\documentclass[letterpaper, 10pt, conference]{ieeeconf}%
\usepackage{amsmath,color}
\usepackage{graphicx}
\usepackage{wrapfig}
\usepackage{amssymb}
\usepackage{mathrsfs}
\usepackage{cancel}
\usepackage{array}
\usepackage{amsmath}
\usepackage{amsfonts}
\usepackage{multirow}
\usepackage{mathtools,amssymb}
\usepackage{algorithmic}
\usepackage{comment}%
\providecommand{\U}[1]{\protect\rule{.1in}{.1in}}
\DeclareMathAlphabet\mathbfcal{OMS}{cmsy}{b}{n}
\IEEEoverridecommandlockouts                              
\begin{document}

\title{Load Balancing in Mobility-on-Demand Systems: \\Reallocation Via Parametric Control Using Concurrent Estimation}
\author{Rebecca M. A. Swaszek$^{1}$ and Christos G. Cassandras$^{1,2}$\\{\small $^{1}$Division of Systems Engineering, $^{2}$Department of Electrical
and Computer Engineering}\\Boston University, Boston, MA 02215, USA\\E-mail:\texttt{\{swaszek,cgc\}@bu.edu \thanks{* Supported in part by NSF under
grants ECCS-1509084, DMS-1664644, CNS-1645681, by AFOSR under grant
FA9550-15-1-0471, by ARPA-E's NEXTCAR program under grant DE-AR0000796 and by
the MathWorks.} }}
\maketitle

\begin{abstract}
Mobility-on-Demand (MoD) systems require load balancing to maintain consistent
service across regions with uneven demand subject to time-varying traffic
conditions. The load-balancing objective is to jointly minimize the fraction
of lost user requests due to vehicle unavailability and the fraction of time
when vehicles drive empty during load balancing operations. In order to bypass
the intractability of a globally optimal solution to this stochastic dynamic
optimization problem, we propose a parametric threshold-based control driven
by the known relative abundance of vehicles available in and en route to each
region. This is still a difficult parametric optimization problem for which
one often resorts to trial-and-error methods where multiple sample paths are
generated through simulation or from actual data under different parameter
settings. In contrast, this paper utilizes concurrent estimation methods to
simultaneously construct multiple sample paths from a single nominal sample
path. The performance of the parametric controller for intermediate size
systems is compared to that of a simpler single-parameter controller, a
state-blind static controller, a policy of no control, and a
theoretically-derived lower bound. Simulation results show the value of state
information in improving performance.

\end{abstract}

%\author{\IEEEauthorblockN{Michael Shell} R. M. A. Swaszek and C. G. Cassandras}
%\section{Introduction}
%\IEEEPARstart{V}{ehicle}

\section{Introduction}

Mobility-on-Demand (MoD) systems, such as those operated by ride-sharing
companies Uber and Lyft, as well as traditional taxi-cab services, are
comprised of a fleet of vehicles, service regions, and users who wish to
travel among or within regions according to dynamic temporal-spatial demand
patterns. Current research promises the advent of fully autonomous vehicles
and Autonomous Taxis (ATs) which are expected to make forays into the MoD
sector \cite{greenblatt2015autonomous}. Service providers of these AT fleets
will face myriad regulatory, safety, and management challenges. The literature
surrounding the fleet management of MoD systems considers how to efficiently
route vehicles among regions to meet current and projected demand, determine
the appropriate fleet size, and guarantee certain levels of system
performance, i.e., average wait times, availability, etc. The load balancing
portion of fleet management can be either proactive in anticipating future
demand or reactive in meeting current demand.

In this paper, we focus on a crucial operational challenge, that of
\emph{proactive load balancing}. This refers to the process of dynamically
redistributing the fleet so as to maintain availability across service areas
to meet future predicted demand.  Clearly, there is a trade-off between
satisfying system-wide customer requests and driving hours logged by empty
vehicles while performing load balancing operations. Whereas profit steers the
driver-centric actions of traditional taxi or Uber drivers, ATs will operate
cohesively under a centralized controller. Rather than many vehicles competing
for users at high-demand areas (e.g., restaurant districts), AT fleet
operators may spread their available vehicles across the broader service area.

Within \emph{proactive} rebalancing there are two general control approaches:
a set of stationary controls based upon historical data and dynamic control
based upon the current state of the system. Reported research results often
demonstrate their control policies using simulation with system parameters
drawn from actual transportation data (taxi, commuter questionnaires)
\cite{zhang2016control},\cite{burghout2015impacts}%
,\cite{brownell2014driverless},\cite{spieser2014toward}%
,\cite{spieser2016vehicle}. A common abstraction is to amalgamate demands from
a region to a single point and utilize queueing structures to represent
available vehicles or waiting customers.

A closed Jackson queueing network model with regards to the vehicles is used
in \cite{zhang2016control} and the proposed controller redistributes ATs by
creating a \textquotedblleft false" user demand rate to force empty ATs from
popular destination regions to popular origin regions to mitigate demand
differences. This control is static in the sense that the \textquotedblleft
false" user rates are set at the beginning of an operational period and do not
change with the state of the system, i.e., the varying availability of ATs
among regions. To address this limitation, a time-driven (as opposed to
event-driven) controller is proposed, which evenly redistributes the AT fleet
among all regions at regular time intervals. Likewise, a MoD system is modeled
in \cite{spieser2014toward} as a closed Jackson network with static
\textquotedblleft false\textquotedblright\ user rates to redistribute vehicles
and explores the trade-off between fleet size and user satisfaction as defined
by average waiting time.

The receding horizon controller introduced in \cite{ramezani2018dynamic} dynamically adjusts false\textquotedblright\ user rates
according to the state of the system.
Greedy state-based controllers are developed in \cite{fagnant2014travel} to
relocate available vehicles by considering the relative abundance of vehicles
and waiting/expected customers in adjacent regions.

Other authors utilize a flow abstraction such that user demand rates induce
vehicle \textquotedblleft flows\textquotedblright\ between regions. In order
to stop some regions from running out while others amass flow, a set of static
controlled rates is proposed in \cite{pavone2012robotic} to redistribute the
flow to ensure stability; \cite{spieser2016vehicle} builds off the
aforementioned work to minimize total system cost in terms of capital costs,
operating costs, and passenger experience, i.e., wait time.

The literature on \emph{reactive} fleet rebalancing focuses primarily on the
vehicle-user assignment problem with an objective to minimize required fleet
size. Whereas in proactive rebalancing, available vehicles are sent empty to
other regions, in reactive rebalancing free vehicles sit idle at their last
destination \cite{alonso2017demand}. Ride-sharing focused papers
\cite{burghout2015impacts},\cite{brownell2014driverless}%
,\cite{vazifeh2018addressing} propose using ATs to operate carpools of varying
sizes and delay times.

In this paper, we use a queueing network model akin to that in
\cite{zhang2016control} in which ATs are discrete entities (as opposed to flow
approximations) with an objective of minimizing the fraction of user requests
that are dropped due to AT unavailability and the fraction of time ATs spend
on load balancing operations driving empty between regions. The optimal
control over an infinite horizon for such a system can be determined using
dynamic programming, but the well-known \textquotedblleft curse of
dimensionality\textquotedblright\ renders such solutions intractable for all
but very small systems.

Our approach is to transform this intractable dynamic optimization problem
into a still difficult but much more manageable parametric optimization
problem where the parameters are \emph{thresholds} in controllers that direct
redistribution based upon the relative quantity of available ATs in and en
route to each region. These thresholds may be tuned to various demand patterns
(rush-hour, high traffic, etc.). We propose two controllers: a time-driven
single-parameter controller and an event-driven multi-parameter controller;
the latter demonstrates superior performance but its parameters require more
effort to tune.

The contribution of the paper is to solve the load balancing problem of a MoD
system by formulating a threshold-based parametric optimization problem and
using concurrent estimation methods \cite{cassandras1999concurrent}%
,\cite{cassandras2009introduction} to estimate the optimal thresholds from a
\emph{single observed sample path} of the queueing network, thus bypassing the
need for repeated trial-and-error. In addition, we derive a lower bound to
assess the performance of our proposed threshold-based control.

%The three major contributions of this paper as follows. First, we solve the load balancing problem of a MoD system by forming a threshold-based parametric optimization problem. Second, we solve for well performing thresholds using a variation of the Standard Clock (SC) Concurrent Estimation (CE) method \cite{cassandras1999concurrent},\cite{cassandras2009introduction}. Finally, we derive a lower bound to assess the performance of our proposed threshold-based control.

%%
Section \ref{System_Model} introduces the MoD system model framework while
Section
%\ref{Optimal_Policy_Simple} describe the challenges of finding the optimal control for simple and general systems, respectively.
\ref{AT_Optimal_Control_Policy} details the challenges associated with finding
an optimal control. Section \ref{Parametric_Controllers} proposes two
parametric threshold-based controllers. 
Section \ref{SC} outlines the Concurrent Estimation techniques of the Standard Clock method as well as the variation used to simulate the AT system.
Section \ref{Section_Lower_Bound}
derives an average best possible performance lower bound and the parametric
controllers' performance is demonstrated in Section \ref{Example} for a
6-region example. Finally, Section \ref{Conclusions} concludes and
highlights potential future work.

%Section \ref{System_Model} introduces the MoD system model framework while Section \ref{Optimal_Policy_Simple} describe the challenges of finding the optimal control for simple and general systems, respectively. Section \ref{Parametric_Controllers} proposes two parametric controllers and Section \ref{SC} gives an overview of concurrent estimation methods. Section \ref{Section_Lower_Bound} derives an average best possible performance lower bound and the parametric controllers' performances are demonstrated in Section \ref{Example} for a 6-region example. Finally, Section \ref{Conclusions} concludes and highlights potential future work.

\section{System Model}

\label{System_Model}

We model a MoD system as a closed Jackson queueing network of $N$ nodes
$\mathcal{N}$ $=\{1,...,N\}$ and $m$ resources representing regions and
vehicles, respectively, similar to the model in \cite{zhang2016control}. We
focus on the load balancing of an urban autonomous taxi fleet and as such
shall refer to vehicles and demands as ATs and
requests, respectively.

\begin{figure}[th]
\centering
\includegraphics[width=1\columnwidth]{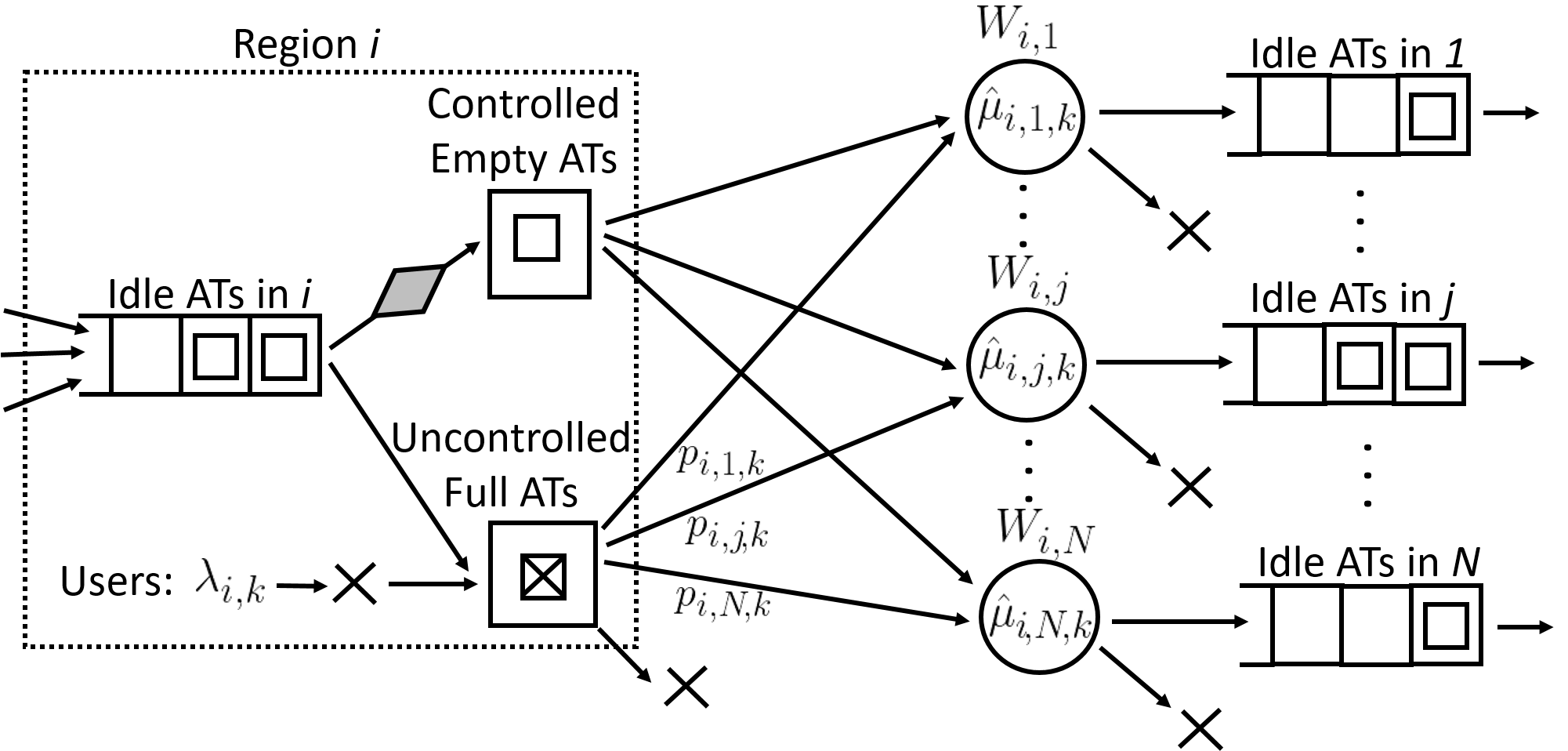}\caption{Region
$i$ consists of users and ATs that are coupled and routed with probability
$p_{i,j,k}$ to a infinite capacity server $W_{i,j}$, then uncoupled and the AT
routed into region $j$. Idle ATs may also be forced by a decision
$\diamondsuit$ to depart empty from the idle queue in $i$ for another
region $j$.}%
\label{region_i}%
\end{figure}

Figure \ref{region_i} shows a region $i$ which consists of a queue
of available ATs. In order to capture time-varying operating demand, routing
and service characteristics, we divide a finite time period $[0,T]$ into $K$
intervals indexed by $k=1,2,...,K$, each of length $I$. Thus, the user request
rate $\lambda_{i,k}$ depends on interval $k=1,2,...,K$. When a user request
occurs, if there is an available idle AT, then the (user, AT) pair, denoted by
$\times$ and $\Box$, respectively, are joined and routed with probability
$p_{i,j,k}$ (including intra-region trips $i=j$) to an infinite-capacity
server $W_{i,j}$ with service rate $\mu_{i,j,k}$. This server captures the
delay experienced by users as they travel from $i$ to $j$. Upon arrival at
region $j$, the pair is separated: the AT is routed with probability 1
to the idle AT queue in region $j$ and the user exits the system. If
there are no available ATs at the time of the user request event, this user
immediately exits the system and incurs a cost. A load balancing controller at
each region, marked by $\diamondsuit$, routes (according to some control
policy) an empty available AT to server $W_{i,j}$ ultimately destined for
region $j$. 

\subsection{State Space}

Let $x_{i}(t)\ \in\{0,1,...,m\}$ be the number of available idle ATs in
region $i \in\mathcal{N}$ and \textbf{\emph{x}}$(t) =[x_{1}%
(t),...,x_{N}(t)]$ be the idle AT vector. Let $y_{i,j}(t)$ be the number of
full ATs (with passengers) en route from $i$ to $j$ and $\mathbf{Y}(t)$ be the
corresponding $N \times N$ matrix populated by $y_{i,j}(t)$. Likewise, let
$z_{i,j}(t)$ be the number of empty ATs in server $W_{i,j}$ and $\mathbf{Z}%
(t)$ be the corresponding $N \times N$ matrix populated by $z_{i,j}(t)$.
Finally let $k(t)\ \in\ \{1,...,K\}$ be the interval that specifies the user
arrival rates $\lambda_{i,k}$, routing probabilities $p_{i,j,k}$, and service
rates $\mu_{i,j,k}$ in effect at time $t$. Thus, the state of the system is
$\mathscr{X}(t)=[\mathbf{x}(t),\mathbf{Y}(t), \ \mathbf{Z}(t),k(t)]$.

\subsection{Events}

\label{events_AT}

The system dynamics are event-driven with the event set $E=E_{U}\cup E_{C}$
where $E_{U}$ and $E_{C}$ contain the uncontrollable and controllable events,
respectively. We define the following uncontrollable event types within
$E_{U}$:

\begin{itemize}
\item $\kappa_{k}$ event: the start of the $k$th interval which prompts a
change of $\lambda_{i,k}$, $p_{i,j,k}$, and $\mu_{i,j,k}$.

\item $\delta_{i,j}$ event: a user request occurs for a trip departing from
region $i$ with destination $j$ (with the possibility that $j=i$). Note
that this request event does not necessitate an AT departure in the case that
region $i$ has no available ATs (see Section \ref{state_dynamics}).

\item $\alpha_{i,j}$ event: a full AT originating from region $i$
arrives at region $j$.

\item $\nu_{i,j}$ event: an empty AT originating from region $i$ arrives
at region $j\neq i$.
\end{itemize}

We also define the following controllable event types within $E_{C}$:

\begin{itemize}
\item $\omega_{i,j}$ event: an empty AT (without passengers) departs from the
idle AT queue in region $i$ destined for region $j\neq i$.

\item $\sigma$ event: a timeout event used for time-driven control.
\end{itemize}

The control policy we select determines when controllable events are
triggered. For example, a controllable $\omega_{i,j}$ event may be triggered
by a timeout or the occurrence of an uncontrollable event resulting in the
state of the system meeting certain criteria.

For a sample path of length $T$ of the MoD system let $\mathbf{e}%
=\{e_{1},...,e_{Q_{T}}\}$ be the observed event sequence, $e_{i}\in E$, with
corresponding event times $\pmb{\tau}=\{\tau_{1},...,\tau_{Q_{T}}\}$ for a
total of $Q_{T}$ events in $[0,T]$. Only at these event times $\tau_{q}$ may
the state of the system change. We may now write the state of the system at
time $\tau_{q}$ as $\mathscr{X}(\tau_{q})\equiv\mathscr{X}_{q}$, with
$x_{i}(\tau_{q})\equiv x_{i,q}$, $y_{i,j}(\tau_{q})\equiv y_{i,j,q}$,
$z_{i,j}(\tau_{q})\equiv x_{i,j,q}$, for all $i \in\mathcal{N}$, where
$q=1,\ldots,Q_{T}$ is the asynchronous event counter.

\subsection{Controls}

A control action in this system consists of forcing an idle AT in some
region $i$ to travel empty to some other region $j\neq i$. This
action depends upon the availability of idle ATs. Let $u_{i,j,q}(x_{i,q})
\in\{0,1,...,x_{i,q}\}$ be the number of empty ATs forced from $i$ to $j$ when
the $q$th event occurs and let $\mathbf{U}_{q}(\mathbf{x}_{q})$ be the
$N\times N$ matrix populated by $u_{i,j,q}(x_{i,q})$. The following
feasibility constraint is required for any control matrix $\mathbf{U}%
_{q}(\mathbf{x}_{q})$:
\begin{equation}
\sum_{j=1}^{N} u_{i,j,q}(x_{i,q})\leq x_{i,q} \quad\forall\ i\in\mathcal{N}%
\end{equation}
For simplicity of notation, let us drop the explicit control dependence on
$x_{i,q}$ and write $u_{i,j,q}$. Note that such controls may be event-driven
(deployed when the state of the system satisfies certain conditions) or
time-driven via the controllable timeout event $\sigma$.

\subsection{State Dynamics}

\label{state_dynamics}

The inventory $x_{i,q}$ of the idle AT queue in region $i$ depends on
both the uncontrollable and controllable events:
\begin{equation}
x_{i,q}=%
\begin{cases}
x_{i,q-1}+1 & e_{q}=\alpha_{j,i}\text{ or }e_{q}=\nu_{j,i}\\
\max\{x_{i,q-1}-1,0\} & e_{q}=\delta_{i,j}\\
x_{i,q-1}-u_{i,j,q} & e_{q}=\omega_{i,j}\\
x_{i,q-1} & \mbox{otherwise}
\end{cases}
\label{x_dyn_AT}%
\end{equation}
where $i,j\ \in\mathcal{N}$. Note that the $\max$ operation prevents the idle
AT queue inventory from falling below 0 in the case that a $\delta_{i,j}$ user
request event occurs and there are no idle ATs, i.e., when a user exits the
system prematurely as in Figure \ref{region_i}.

The number of full ATs $y_{i,j,q}$ en route from region $i$ to
region $j$ evolves according to:
\begin{equation}
y_{i,j,q}=%
\begin{cases}
y_{i,j,q-1}+1 & e_{q} = \delta_{i,j} \mbox{ and } x_{i,q-1} > 0\\
y_{i,j,q-1} - 1 & e_{q}= \alpha_{i,j}\\
y_{i,j,q-1} & \mbox{otherwise}\\
\end{cases}
\label{y_dyn_AT}%
\end{equation}

Likewise the number of empty ATs en route from $i$ to $j$ evolves according
to:
\begin{equation}
z_{i,j,q}=%
\begin{cases}
z_{i,j,q-1} + u_{i,j,q} & e_{q} = \omega_{i,j}\\
z_{i,j,q-1} - 1 & e_{q}= \nu_{i,j}\\
z_{i,j,q-1} & \mbox{otherwise}\\
\end{cases}
\label{z_dyn_AT}%
\end{equation}

At time instants $t=k\mathcal{I}$, $k\in\{1,\ldots,K\}$ the interval index $k$
changes upon the occurrence of a $\kappa_{k}$ event, thus causing the
time-varying parameters to change for all regions: $k_{q}=k_{q-1}+1$

\subsection{Objective Function}

The primary concern in the relevant literature of MoD systems is to meet all
or most of the user requests while also minimizing load balancing operations.
We formalize this trade-off in an objective function which minimizes a
weighted sum of the probability that a user's request does not result in an AT
departure (i.e., is rejected as in Figure \ref{region_i}) and the
probability that an AT is driving empty performing load balancing operations.
The evaluation of these probabilities is generally infeasible because of the
fast growth of the state space, rendering this task intracatable except for
the simplest of systems. In order to assess the effect of controls $u_{i,j,q}$
on system performance, we replace these probabilities with estimates
consisting of the fraction of rejected requests and total time spent driving
empty in a sample path over $[0,T]$. Let $\pmb{\rho}_{T}^{i}=\{\rho
_{i,1},...,\rho_{i,D_{T}^{i}}\}$ be the event times of all $\delta
_{i,j},\ j\in\mathcal{N}$ user request events at region $i$ where
$D_{T}^{i}$ is the total observed number of such events, and let
$\mathbf{1}[\cdot]$ be the usual indicator function. As defined in Section
\ref{events_AT}, there are a total of $Q_{T}$ events in $[0,T]$ and the $q$th
event occurs at time $\tau_{q}$.

The objective function we define is:
\begin{multline}
J(\!\mathscr{X}_{0})=E\bigg [w\frac{\sum_{i=1}^{N}\sum_{q=1}^{D_{T}^{i}%
}\mathbf{1}[x_{i}(\rho_{i,q})=0]}{\sum_{i=1}^{N}D_{T}^{i}}\\
+(1-w)\frac{\sum_{q=1}^{Q_{T}}\sum_{i=1}^{N}\sum_{j=1}^{N}(\tau_{q}-\tau
_{q-1})z_{i,j,q}}{Tm}\!\bigg ] \label{J_intro}%
\end{multline}
where $w\in(0,1]$ is a weight coefficient. The first part of (\ref{J_intro})
refers to the fraction of rejected users unable to obtain an available idle AT
(the numerator is the total number of all user request events). The second
part refers to the fraction of vehicle-hours when ATs drive empty (the
denominator is the total vehicle-hours driven by the $m$ vehicles in the fleet
over $[0,T]$). The weight coefficient $w$ is used to quantify the trade-off
between customer satisfaction and load balancing effort. We exclude $w=0$ as
the optimal control is trivial when customer satisfaction is irrelevant. If
$w=1$, the optimal control is still difficult to determine: although load
balancing may not be a direct cost, the unavailability of ATs while performing
load balancing operations creates an indirect cost. Note that the objective
function in (\ref{J_intro}) is properly normalized so as to give values
bounded by $[0,1]$ and the weight $w$ creates a convex combination of the two
objective components.

The optimization problem we formulate based on (\ref{J_intro}) is to determine
a control policy $\mathbf{U}(\mathscr{X})$ so as to minimize this objective:
\begin{equation}
J^{\ast}(\mathscr{X}_{0})=\min_{\mathbf{U}(\mathscr{X})}J(\!\mathscr{X}_{0})
\end{equation}

\section{Optimal Control Policy} \label{AT_Optimal_Control_Policy}

Let us first consider a simpler version of the problem where we assume that
user arrivals occur according to a Poisson process with fixed rate
$\lambda_{i}$ and that each infinite-capacity server $W_{i,j}$ has
exponentially distributed service times with mean service time $\frac{1}%
{\mu_{i,j}}$. Thus, the MoD system is described by a finite-dimensional
continuous time Markov chain. As an example, the simplest possible such system
corresponds to $N$=2 and $m$=1, in which case the single AT is in one of eight
possible states: idle in either region 1 or 2, en route with a passenger
in one of the four infinite capacity queues $W_{1,1},W_{1,2},W_{2,1},W_{2,2}$
or en route empty performing load balancing in $W_{1,2}^{E},W_{2,1}^{E}$ (with
the subscript $E$ denoting \textquotedblleft empty\textquotedblright). For
this simple system, the optimal control policy may be found analytically as
shown in  Appendix A and B. % in \textbf{[arXiv LINK HERE]}. 
For any larger
system, we turn to dynamic programming (DP) \cite{bertsekas2005dynamic} to
determine the optimal control policy for the average cost over the infinite
horizon. The cost $C(\mathscr{X})$ associated with being in state
$\mathscr{X}=[\mathbf{x},\mathbf{Y},\ \mathbf{Z}]$ is as follows:
\begin{multline}
C(\mathscr{X})=w\sum_{i=1}^N \frac{\sum_{j=1}^N \lambda_{i,j}\mathbf{1}[x_{i}=0]}%
{\sum_{i=1}^N \sum_{j=1}^N (\lambda_{i,j}+\mu_{i,j})}\\
+(1-w)\frac{\sum_{i=1}^N \sum_{j=1}^N %
z_{i,j}}{m}\label{state_x_cost}%
\end{multline}
Let $S$ and $U$ denote the state and control spaces, respectively, let
$\bar{g}(i,u)$ be the cost of state $i$ as defined by \eqref{state_x_cost},
and suppose that the state space cardinality is $|S|=n$. Let $P_{i,j}(u)$ be
the transition probability from state $i$ to state $j$ under control $u$. We
wish to solve for the optimal cost function $J^{\ast}$ which satisfies for all
$i\in S$ the deterministic form of Bellman's equation
\cite{bertsekas2005dynamic}:
\begin{equation}
J^{\ast}(i)=\min_{u\in U(i)}\bar{g}(i,u)+\sum_{j=1}^{n}P_{i,j}(u)J^{\ast
}(j)\label{Bellman}%
\end{equation}

Linear programming or policy iteration may be used to solve \eqref{Bellman}.
In the linear programming formulation, we use the differential cost vector
$[h(1),...h(n)]$ (necessary in the absence of discounting) and the optimal
cost $\bar{J}$ as the decision variables \cite{bertsekas2005dynamic}:
\begin{multline}
\max_{\bar{J},h(1),...h(n)}\bar{J} \qquad \qquad \mbox{s. t. } \\
\bar{J}+h(i)\leq\bar{g}(i,u)+\sum_{j=1}^{n}P_{i,j}(u)h(j)\qquad\forall\ i\in
S,\ u\in U(i)
\end{multline}
The resulting optimal steady-state distribution is given by the Lagrange
multipliers of $[h(1),...,h(n)]$ and the binding constraints indicate optimal
policies. As an example, the full linear programming problem for the $N=2$,
$m=1$ system may be found in Appendix C.% of \textbf{[arXiv LINK HERE]}.

Likewise, the policy iteration method solves a linear system of equations
$\mathbf{h}_{\mu}+\mathbf{e}\bar{J}_{\mu}=\mathbf{\bar{g}}_{\mu}%
+\mathbf{P}_{\mu}\mathbf{h}_{\mu}$ in which $\mathbf{h}_{\mu}$ and
$\mathbf{\bar{g}}_{\mu}$ are the $n\times1$ vectors of $\bar{g}(i)$ and $h(i)$
under some policy $\mu$, $\mathbf{P}_{\mu}$ is the transition probability
matrix populated by $P_{i,j}(\mu)$, and $\mathbf{e}$ is a $n\times1$ vector
with $e_{i}=1$ for all $i=1,\ldots,n$. Policy iteration is the more efficient
method to find the optimal control as it solves smaller linear programs and
requires less memory, although convergence to the optimal policy may be slow,
depending on the initial policies chosen.

However, both methods are limited to small systems as the state space grows
combinatorially with the number of regions and ATs. For a system with
$N$ regions and $m$ ATs, the cardinality of the state space is
$\binom{m+2N^{2}-1}{2N^{2}-1}$ and the cardinality of states with at least
$N-1$ controls is: $\binom{m+2N^{2}-1}{2N^{2}-1}-\binom{m+2N^{2}-N-1}%
{2N^{2}-N-1}$. For example, an AT system with $N$=6 regions and $m$=50
taxis has over 3$\times10^{34}$ states. Utilizing sparse matrices and shared cloud computing
facilities, the largest possible systems that we have been able to analyze using
the policy iteration method consists of $N$=2 regions and $m$=10 ATs and
of $N$=3 regions and $m$=5 ATs. The latter required only 31 GB of
memory, while solving for $N$=3 regions and $m$=6 ATs required more than
the 256 GB allotted to a full node of 28 cores. This provides the motivation
for seeking alternative control policies and assessing their performance using
the fraction estimates in \eqref{J_intro}.

\section{Parametric Control Policies}

\label{Parametric_Controllers}

In this section, we introduce a parametric controller for larger MoDs
expanding upon the framework of the \textquotedblleft real-time" controller in
\cite{zhang2016control}; this \textquotedblleft real-time\textquotedblright%
\ controller rebalances ATs evenly among regions every half hour using an
integer linear programming approach to minimize the expected load balancing
time. Unlike our loss model in which user requests may be rejected if there
are no available ATs, the model in \cite{zhang2016control} includes a queue
for waiting users such that the number of ATs associated with a region is the
sum of available ATs and ATs en route with the number of users queued up
waiting for an AT subtracted from the latter.

Let $\Theta=[\theta_{1},...,\theta_{N}]$ be a parameter vector with
$\theta_{i}\geq0$, $i=1,\ldots,N$ and $\sum_{i=1}^{N}\theta_{i}\leq m$. This
vector defines a \textquotedblleft fill to\textquotedblright\ level for each
of the $N$ regions; this is akin to $(s,S)$ threshold policies in supply
chain and inventory management \cite{bertsekas2005dynamic} where $s$ is a
\textquotedblleft fill to\textquotedblright\ level such that when an inventory
drops below it, a supply request is triggered (similarly, crossing $S$ from
below triggers a request to stop the supply process). Note that each interval
$k=1,2,...,K$ may have its own set of parameters to account for different
request patterns and traffic conditions. For simplicity of notation, let
$a_{i}(t)$ be the total number of ATs available at or en route to region
$i$ at time $t$:
\begin{equation}
a_{i}(t)=x_{i}(t)+\sum_{j=1}^{N}\big (y_{j,i}(t)+z_{j,i}(t)\big )
\end{equation}

Furthermore let us define a quantity $D_{i}(t)$ that is the \emph{supply} of
available excess ATs if positive or the \emph{demand} for ATs if negative in
region $i:\ \ \ D_{i}(t)=\min\{a_{i}(t)-\theta_{i},x_{i}(t)\}$ Note that
$D_{i}(t)$ is an integer quantity as $\theta_{i},a_{i}(t),$ and $x_{i}(t)$ are
all integers. In order for feasible AT redistribution actions to be triggered,
it is a necessary condition that the overall supply must exceed the demand in
the following inequality:
\begin{equation}
\sum_{k\in\{i\in\mathcal{N}|D_{i}(t)>0\}}\mkern-18mu\mkern-18muD_{k}%
(t)\qquad\geq\sum_{j\in\{i\in\mathcal{N}|D_{i}(t)\leq0\}}%
\mkern-18mu\mkern-18muD_{j}(t)\label{need_have}%
\end{equation}
%\begin{equation}
%\sum_{k\in\{i\in\mathcal{N}|a_{i}(t)>\theta_{i}\}}\mkern-18mu\mkern-18mu\min
%\{a_{k}(t)-\theta_{k},x_{k}(t)\}\quad\geq\mkern-18mu\sum_{j\in\{i\in
%\mathcal{N}|a_{i}(t)\leq\theta_{i}\}}\mkern-18mu\mkern-18mu\theta_{j}%
%-a_{j}(t)\label{need_have}%
%\end{equation}
This simply asserts that there is an adequate number of available ATs in
regions which are above their \textquotedblleft fill-to\textquotedblright%
\ levels specified in $\Theta$ which can be used to supply those regions
whose queues of available ATs are below their \textquotedblleft
fill-to\textquotedblright levels.

Assuming for the moment that there exists a well-defined mechanism for
triggering a process to redistribute ATs among regions (further
discussed in Sections \ref{AT_1} and \ref{AT_N_plus_1}), this process consists
of the following integer linear program with decision variables $u_{i,j}%
\in\{0,1,2,\ldots\}$:
\begin{equation}
\min_{u_{i,j,}\text{ }i,j\in\mathcal{N}}\sum_{i=1}^{N}\sum_{j=1}^{N}%
\frac{u_{i,j}}{\mu_{i,j}} \label{MINLP_event}%
\end{equation}%
\[
\text{s. t. \ }%
\begin{aligned}[t] \theta_i \le a_i(t) + \sum_{j=1}^{N} \big ( u_{j,i} - u_{i,j} \big ) & & i & \in \mathcal{N} \\
\sum_{j=1}^{N} u_{i,j} \le x_{i}(t) & & i & \in \mathcal{N} \\
u_{i,j} \in \{0,1,2,... \} && \  i,j& \in \mathcal{N}
\end{aligned}
\]

The objective function of (\ref{MINLP_event}) minimizes empty vehicle driving
time; the first constraint requires the intended inventory for each
region to meet or exceed \textquotedblleft fill to\textquotedblright%
\ levels $\theta_{1},...,\theta_{N}$, and the second constraint maintains feasibility.

In order to bypass the difficulty of integer programming, we rewrite
\eqref{MINLP_event} as a relaxed linear program in the form of a minimum cost
flow problem in which regions indexed by $i$ with positive \emph{supply}
$D_{i}(t)$ are sources and negative \emph{demand} $D_{i}(t)$ are sinks.
\begin{equation}
\min_{u_{i,j,}\text{ }i,j\in\mathcal{N}}\sum_{i=1}^{N}\sum_{j=1}^{N}%
\frac{u_{i,j}}{\mu_{i,j}}\label{min_cost_flow}%
\end{equation}%
\[
\text{s. t. \ }%
\begin{aligned}[t] D_i(t) \ge \sum_{j=1}^{N} \big ( u_{i,j} - u_{j,i} \big ) & & i & \in \mathcal{N} \\
u_{i,j} \ge 0 && \  i,j& \in \mathcal{N}
\end{aligned}
\]
Note that the single constraint in \eqref{min_cost_flow} encompasses both
constraints in \eqref{MINLP_event}: $(i)$, if $D_{i}(t)=a_{i}(t)-\theta_{i}%
$, the constraint in \eqref{min_cost_flow} is identical to the first
constraint of \eqref{MINLP_event}; $(ii)$ if $D_{i}(t)=x_{i}(t)$, then $i$ is
a source such that no flow is directed to it, i.e., $\sum_{j=1}^{N}u_{j,i}=0$,
therefore, the constraint in \eqref{min_cost_flow} becomes identical to the
second constraint of \eqref{MINLP_event}. We recover the integer solution to
\eqref{MINLP_event} from this linear program since integer solutions are a
property of minimum cost flow linear programs with integer sink and source
quantities \cite{BazaraaM.S2010Lpan}.

A potential problem associated with this policy is that, depending upon the
frequency at which these control actions are implemented, there could be empty
ATs cycling inefficiently such that $z_{i,j}(t)>0$ and $z_{j,i}(t)>0$
simultaneously. To address this issue, we introduce an additional parameter
$\Omega$ which acts as either a scalar \textquotedblleft
timeout\textquotedblright\ in a time-driven controller in Section \ref{AT_1}
or as an integer threshold in an event-driven controller in Section
\ref{AT_N_plus_1}. In either case, the role of $\Omega$ is to define the
mechanism that triggers the process of solving (\ref{min_cost_flow}) so as to
determine the values of the control variables $u_{i,j}$, $i,j\in\mathcal{N}$.

\subsection{Single Scalar Parameter Time-Driven Controller}

\label{AT_1}

For a time-driven controller, let the scalar parameter $\Omega\in(0,\infty)$
be associated with the timeout event $\sigma$ we defined in Section
\ref{events_AT} such that event $\sigma$ triggers the control actions in
\eqref{min_cost_flow} every $\Omega$ time units. The \textquotedblleft
real-time" controller proposed in \cite{zhang2016control} can be recovered
from our controller by setting all $\theta_{i}=\lfloor\frac{m}{N}%
\rfloor\ ,\ i\in\mathcal{N}$ in \eqref{MINLP_event} and $\Omega=30$ minutes.
This time-driven single scalar parameter controller is quite effective, but
does not account for the system-specific demand rates. For example consider a
heterogeneous two-region system -- this control of distributing the fleet
equally between them would perform badly if one region experienced far greater
requests than the other.

\subsection{$N+1$ Integer Parameter Event-Driven Controller}

\label{AT_N_plus_1}

In order to avoid the aforementioned inefficiencies associated with equally
spreading out the fleet among regions of dissimilar demand, we define an
\emph{event-driven} controller to trigger a solution of problem
(\ref{min_cost_flow}) whenever $a_{i}(t)$, the number of ATs in or en route to
a region, drops below the threshold $\theta_{i}$ by some amount. In the
simplest case, this occurs as soon as $a_{i}(t)<\theta_{i}$, corresponding to
a greedy mechanism that pulls a single empty AT from the nearest region $j$
with $a_{j}(t)>\theta_{j}$ and setting $u_{i,j}=1$. However, since we have at
our disposal a central controller with full information of all region
states and AT locations, we can do better than that as explained next.

Let us redefine $\Omega\in\{1,2,..,m\}$ as an integer-valued threshold
parameter used to trigger the control actions resulting from a solution of
(\ref{min_cost_flow}) whenever the condition $\theta_{i}-a_{i}(t)>\Omega$ is
satisfied. This will effectively send a total of $\theta_{i}-a_{i}(t)$ empty
ATs to region $i$. However, this is still a region-specific
control that may be inefficient over the system as a whole. As such, we
instead trigger the control when the sum-positive difference between the
\textquotedblleft fill-to\textquotedblright\ levels $\theta_{i}$ and available
or en route ATs $a_{i}(t)$ surpasses threshold $\Omega$ \emph{across all
regions} for a global centralized control policy triggered by:
\begin{equation}
\sum_{j\in\{i\in\mathcal{N}|a_{i}(t)<\theta_{i}\}}%
\mkern-18mu\mkern-18mu\big (\theta_{j}-a_{j}(t)\big )\quad>\quad
\Omega\label{Omega}%
\end{equation}

\begin{figure}[th]
\centering
\includegraphics[width=7 cm]{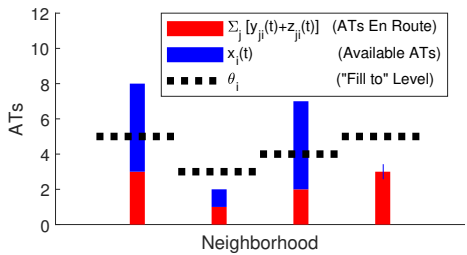}\caption{Control parameters
$\theta_{i}$ represent a \textquotedblleft fill to" level for the number of
ATs available at or en route to region $i$. Below this threshold
represents need; above this threshold represents excess inventory to possibly
send elsewhere.}%
\label{nplus1_statespace}%
\end{figure}

For example, consider the state of a system with $N=4,$ $m=20$ in Fig.
\ref{nplus1_statespace}. The control parameter vector $\Theta=[5,3,4,5]$ is
hashed in black and let $\Omega=2$. After each event, the controller checks
inequalities \eqref{need_have} and \eqref{Omega}; if both hold, then
$\omega_{i,j}$ events are induced as per \eqref{min_cost_flow}. In this
example, \eqref{need_have} holds with (0+1+0+2) $\leq$ (3+0+3+0) and
\eqref{Omega} holds with (0+1+0+2) $>$ 2. 

This controller is triggered when one of the events defined in Section
\ref{events_AT} causes a change in the value of some $a_{i}(t)$ in (\ref{Omega})
so that the sum crosses $\Omega$ either from below or from above. This may
happen in two ways: $(i)$ When event $\delta_{i,j}$ occurs, i.e., a user
requests to go from region $i$ to $j$, which alters $x_{i}(t)$ and
$y_{i,j}(t)$ as per the state dynamics \eqref{x_dyn_AT} and \eqref{y_dyn_AT},
respectively (note that this assumes \eqref{need_have} already holds). $(ii)$
Event $\alpha_{j,i}$ or $\nu_{j,i}$ occurs, i.e., a full (respectively, empty)
AT arrives at region $i$ and alters the value of $x_{i}(t)$ and
$y_{j,i}(t)$ (respectively, $z_{j,i}(t)$).

While the time-driven controller triggers control actions after a
predetermined length of time regardless of the state of the system, the
event-driven controller is only triggered when the inventory levels fall
sufficiently low across all regions.

\subsection{$N^{2}$ Parameter Static Controller}

Both of the previously described parametric controllers require tuning the
threshold parameters via simulation or through a data-driven on-line
adaptation process. We shall compare in Section \ref{Example} these
controllers to an alternative simpler parametric controller introduced in
\cite{pavone2012robotic} whose parameters are determined by linear programming
and rely solely on the model parameters $\lambda_{i,j}$ and $\mu_{i,j}$. This
time-invariant and state-blind control sends empty ATs from $i$ to $j$ at a
static rate $r_{i,j}$, $i,j\ \in\mathcal{N}$ such that there are $N^{2}$ rate
parameters. These static rate parameters are determined by the following
linear program that minimizes empty travel time and seeks to equal the inflow
and outflow of ATs at each region:
\begin{equation}
\min_{r_{i,j}}\sum_{i=1}^{N}\sum_{j=1}^{N}\frac{r_{i,j}}{\mu_{i,j}}
\label{MINLP}%
\end{equation}%
\[
\text{s. t. \ }%
\begin{aligned}[t] \sum_{j=1}^{N}(\lambda_{i,j} + r_{i,j}) = & \sum_{j=1}^{N} (\lambda_{j,i} + r_{j,i} ) & i & \in \mathcal{N} \\ r_{i,j} \ge & \ \ 0 & i,j & \in \mathcal{N} \end{aligned}
\]

This linear program is always feasible. This is easily seen by taking, for
example, $r_{i,j}=\lambda_{j,i}\ \forall\ i,j\in\mathcal{N}$, which satisfies
the first constraint by sending back empty taxis at the same rate and
satisfies the second constraint as $\lambda_{i,j}\geq0$. This control is
state-blind since it does not depend on the current state or even the fleet
size. Note that in the case that such a controller as a comparison baseline at
time $t$ directs an AT to leave region $i$ for some $j$ but $x_{i}(t)=0$
(i.e., region $i$ does not have any idle ATs) neither a penalty nor an event
$\omega_{i,j}$ occur. We use this state-blind controller to demonstrate the
importance of a dynamic (state-driven) control strategy in the example system
of Section \ref{Example}.

\section{Concurrent Estimation Methods} \label{SC}

 In order to find well performing parameters for both the time-driven and event-driven controllers, we turn to concurrent estimation methods to simulate the effects of many different sets of parametric controls from a single sample path and use the resulting performances to direct the evolution of the next set of control parameters to test.

 Concurrent estimation/simulation methods help in problems such as this AT fleet load balancing problem where we wish to test the performance of discrete parameters $\pmb{\theta}$ whose effects on the objective function $J(\pmb{\theta})$ are difficult to determine analytically. The off-line brute force approach to simulate each of the possible sets of parameters is very inefficient and time consuming. On-line each possible parameter set could be tested in a trial and error process, but this too is inefficient and may be impossible given the available time. In the time-driven control, let $\pmb{\theta} = [\Omega]$ and in the event-driven control let $\pmb{\theta} = [\theta_{1},...,\theta_{N},\Omega]$. 

 The off-line Standard Clock (SC) method of concurrent simulation constructs sample paths from a nominal simulated sample path in order to simulate the performance $J(\pmb{\theta})$ under various policies $\pmb{\theta}$ with the drawback of possibly simulating many fictitious events. Standard Clock is used to simulate stochastic timed automata with a Poisson clock structure such as our AT system which is an event-driven system with exponentially distributed event lifetimes. In each state $x$ there is a set of feasible events $\Gamma(x)$. Once the system is in state $x$ it stays there for a exponential amount of time with mean $\frac{1}{\Lambda(x)}$ where $\Lambda(x) = \sum_{i \in \Gamma(x)} \lambda_i$ and $\lambda_i$ is the rate at which event $i$ occurs. The distribution of the next triggering event $i$ is given by 
 \begin{equation}
 p(i,x) = \frac{\lambda_i}{\Lambda(x)} \ , i \in \Gamma(x)
 \label{pi}
 \end{equation}
One may uniformize the Markov chain by introducing a uniform rate: 
 \begin{equation}
 \gamma \ge \Lambda(x)
 \label{unirate}
 \end{equation}
 and replace $\Lambda(x)$ with this rate $\gamma$ for all states $x$ and force the additional probability flow $[\gamma - \Lambda(x)]$ to be a ``fictitious" event that is a self loop back to $x$ leaving the state unchanged. 

 The Standard Clock method builds off this uniformized model and chooses $\gamma = \Lambda$ with 
 \begin{equation}
 \Lambda = \sum_{i \in \mathcal{E}} \lambda_i
 \end{equation} 
 where $\mathcal{E}$ is the set of all events. Note that this choice of $\gamma$ has the potential to be much larger than any $\Lambda(x)$ thus forcing the additional probability flow $[\Lambda - \Lambda(x)]$ to also be large for many states $x$. This causes many fictitious events which slow down simulations. With this new uniform rate $\gamma$ the triggering event probability is now \begin{equation}
 p_i = \lambda_i/\Lambda
 \end{equation}
 which could potentially be much smaller than \eqref{pi}. The occurrence of the $k$th event $E_k$ in the event set $\mathcal{E} = {1, . . . , N}$ may be determined by a random number $U_k$:
 \begin{equation}
 E_k=%
 \begin{cases}
 1& \mbox{if } 0 \le U_k \le \lambda_1/\Lambda  \\
 2 & \mbox{if } \lambda_1/\Lambda < U_k \le (\lambda_1+\lambda_2)/\Lambda  \\
 \cdots & \\
 N & \mbox{if } (\lambda_1+...+\lambda_{N-1})/\Lambda < U_k \le 1  \\ 
 \end{cases}
 \label{11.135}
 \end{equation}
 By construction \eqref{11.135} allows all events to occur in \textit{any} state; however, if the event $E_k$ is not feasible at state $x$ than it is a fictitious event i.e. a self-loop. 

 The following steps outline the method of constructing $M+1$ sample paths concurrently where $f_m(\cdot)$ is the state transition function \cite{cassandras2009introduction}:

 \fbox{\begin{minipage}{22em}

 1. CONSTRUCT STANDARD CLOCK:

 $\quad$ $\{V_1, V_2, . . .\}, V_k  \sim 1 -e^{-t},t > 0$

 $\quad$ For every constructed sample path $m = 0, 1, . . . , M$:

 2. DETERMINE TRIGGERING EVENT $E_m$ BY (\ref{11.135})

 3. UPDATE STATE $X_m:$ 

 $\qquad X_m := f_m(X_m, E_m)$

 4. RESCALE INTEREVENT TIME $V$: 

 $\qquad V_m = V /\Lambda_m$

 \end{minipage}}

 \vspace*{.1in}
 \noindent
 Note that in Step 1 Standard Clock generates inter-event times from an exponential distribution with parameter 1, thus for models with event rate $\Lambda \ne 1$ Step 3 adjusts the sequence $\{V_1, V_2, . . .\}$ by rescaling $V_k(\Lambda) = V_k/\Lambda$.

 \subsection{Variation of Standard Clock}

 We introduce a variation on the Standard Clock method which takes advantage of models in which there is a substantial subset of events that are feasible across all states $x$ and creates fewer fictitious events with a drawback of slightly more calculations on per event. 

 Let $\xi \subseteq \mathcal{E}$ be the subset of the $C$ events that are feasible in every state $x$: 
 \begin{equation}
 \xi = \{ i \in \mathcal{E} | i \in \Gamma(x) , \forall \ x \}
 \label{xi}
 \end{equation}
 and let $\xi_x$ be the other feasible events in state $x$:
 \begin{equation}
 \xi_x = \{ i \in \mathcal{E} | i \in \Gamma(x) , i \notin \xi \}
 \label{xi_x}
 \end{equation}
 Further, set the uniformized rate $\gamma$ as in \eqref{unirate}  
 \begin{equation} \label{SC_max_rate}
 \gamma = \max_x \Lambda(x)
 \end{equation}
 Similar to \eqref{11.135}, for state $x$ the range of $U_k$ is 
 partitioned to the $C$ common events and the remaining events in $\xi_x$
 (plus fictitious events if there is any remaining range for $U_k$). 
 Let $B$ be the probability that the event is within this common event set $\xi$:
 \begin{equation}
 B = \frac{\sum_{i \in \xi}\lambda_i}{\gamma}
 \label{B}
 \end{equation}
 and let $B_x$ be the probability that the event is feasible in state $x$:
 \begin{equation}
 B_x = \frac{\sum_{i \in \Gamma(x)}\lambda_i}{\gamma}
 \label{B_x}
 \end{equation}
 (Clearly at least one $B_x=1$; states with $B_x<1$ have the potential for fictitious events to occur.)

 If $U_k \leq B$, then the event is in $\xi$
 \begin{equation}
 E_k=%
 \begin{cases}
 1& \mbox{if } 0 \le U_k \le \lambda_1/\gamma  \\
 2 & \mbox{if } \lambda_1/\Lambda < U_k \le (\lambda_1+\lambda_2)/\gamma  \\
 \cdots & \\
 C & \mbox{if } (\lambda_1 \! + \! ... \! + \! \lambda_{C-1})/\gamma  \! < U_k  \! \le  B  \\ 
 \end{cases}
 \label{Uk_small}
 \end{equation}
 If $U_k > B$, then event $E_k$ is either in $\xi_x$ or fictitious
 \begin{equation}
 E_k=%
 \begin{cases}
 \xi_x(1) & \mbox{if } B < \! U_k \! \le \! B \! + \! \frac{\lambda_{\xi_x(1)}}{\gamma}  \\
 \cdots & \\
 \xi_x (\big | \xi_x \big | ) & \mbox{if }  B \! + \! \frac{\lambda_{\xi_x(1)}+...+\lambda_{\xi_x ( | \xi_x  | -1)}}{\gamma} \! < \! U_k \! \le B_x \\ 
 \mbox{Fictitious} & B_x < U_k
 \end{cases}
 \label{Uk_big}
 \end{equation}

 Let \eqref{xi},\eqref{B}, and \eqref{Uk_small} be pre-calculated as the common events in $\xi$ will occur often while let \eqref{xi_x}, \eqref{B_x}, and \eqref{Uk_big} be calculated on a need basis as the large state space renders them impossible to save for all states $x$. 

 The following four steps construct $M+1$ concurrent sample paths:

 \fbox{\begin{minipage}{22em}

 1. CONSTRUCT STANDARD CLOCK:

 $\quad$ $\{V_1, V_2, . . .\}, V_k  \sim 1 -e^{-t},t > 0$

 $\quad$ For every constructed sample path m $= 0, 1, . . . , M$:

 2. DETERMINE TRIGGERING EVENT $E_m$  

 \quad If $U_k \le B$ \ \ BY: \eqref{Uk_small} 

 \quad Else \qquad \ \ \ \ BY: \eqref{Uk_big}

 3. UPDATE STATE $X_m:$ 

 $\quad X_m := f_m(X_m, E_m)$

 4. RESCALE INTEREVENT TIME $V$: 

 $\quad V_m = V /\Lambda_m$

 \end{minipage}}

 \vspace*{.1in}
 \noindent
 The advantage of this variation is that by setting a smaller $\gamma$ many fewer fictitious events occur at the cost of the additional calculations in \eqref{Uk_big}. This method works well for our AT model as the set of common events -- all user request events -- is quite large: $\xi = \delta_{k,j}, \ k,j \in \mathcal{N}$ such that the probability $B$ of an event being within the common event set is:  
 \begin{equation}
 B = \frac{\sum_{i,j}\lambda_{i,j}}{\gamma}
 \end{equation}
with $\gamma$ set as per \eqref{SC_max_rate} to be: 
 \begin{equation}
 \gamma = \sum_{i,j} \lambda_{i,j} + m \max_{i,j} \mu_{i,j}
 \end{equation}

 Table \ref{Infeasible_Events} shows the minimum and maximum percent of fictitious events possible for the Standard Clock method and variation thereof (see the example in Section \ref{Example}). 

 {\renewcommand{\arraystretch}{1.7}
 \begin{table}[h]
 \centering
 \caption{Percent Fictitious Events}
 \begin{tabular}{c|c|c|}
 \cline{2-3}
                                      & Minimum & Maximum \\ \hline
 \multicolumn{1}{|l|}{SC} &    $ \! \! \! 1 \text{-}  \frac{\sum_{i,j} \lambda_{i,j} + m \max_{i,j} \mu_{i,j} }{\sum_{i,j}(\lambda_{i,j} + m \mu_{i,j})}$    &   $\! \! \! 1 \text{-} \frac{\sum_{i,j} \lambda_{i,j}}{\sum_{i,j}(\lambda_{i,j} + m \mu_{i,j})}$     \\ \hline
 \multicolumn{1}{|l|}{Variation}      &    0     &     $\! \! \! 1 \text{-} \frac{\sum_{i,j}\lambda_{i,j}}{\sum_{i,j} \lambda_{i,j} + m \max_{i,j} \mu_{i,j} }$     \\ \hline
 \end{tabular}
 \label{Infeasible_Events}
 \end{table}}

We use concurrent simulation to find well performing control parameter vectors via a greedy iterative search and a broader random search of the control space. 

\subsection{Greedy Iterative Search} \label{Greedy_Iterative_Search_AT}

 In an iterative process, we find well performing control vectors by following a greedy process of following vector changes that lead to better performances. For example, consider the time-driven controller with its single parameter $\Omega$ starting at some value $\Omega^0$. We may use concurrent simulation to compare how the slight control deviations such as $\Omega^0-1$, and $\Omega^0+1$ perform and direct the next controller: $\Omega^1 = \mbox{argmin}_{\Omega} \{ J(\Omega^0-1),J(\Omega^0),J(\Omega^0+1) \}$. Iterated again and again with $\Omega^{k+1} = \mbox{argmin}_{\Omega} \{ J(\Omega^k-1),J(\Omega^0),J(\Omega^0+k) \}$ and with an increasing simulation length $T$, this iterative method gravitates towards well performing controls. As this requires simulation, it cannot find the optimal control, but hones in on a flat region of the objective function. 
 
 Likewise for the event-driven controller with $N$+1 control parameters, the parameter vector may start as $[\theta_1^0,\theta_2^0,...\theta_N^0,\Omega^0]$, and the first iteration concurrently creates sample paths for the following control parameters: 
$$[\theta_1^0+1,\theta_2^0,...\theta_N^0,\Omega^0]$$
$$[\theta_1^0,\theta_2^0+1,...\theta_N^0,\Omega^0]$$
$$...$$
$$[\theta_1^0,\theta_2^0,...\theta_N^0+1,\Omega^0]$$
$$[\theta_1^0,\theta_2^0,...\theta_N^0,\Omega^0+1]$$
$$[\theta_1^0-1,\theta_2^0,...\theta_N^0,\Omega^0]$$
$$[\theta_1^0,\theta_2^0-1,...\theta_N^0,\Omega^0]$$
$$...$$
$$[\theta_1^0,\theta_2^0,...\theta_N^0-1,\Omega^0]$$
$$[\theta_1^0,\theta_2^0,...\theta_N^0,\Omega^0-1]$$
At the end of the iteration, the next control parameter picked as the next nominal path would be that which performed the best such that each iterate follows the path of steepest descent.

\subsection{Random Search of the Control Space} \label{Random_Search_AT}

The greedy iterative process as described in the previous subsection searches locally within a small section of the control space and the final selected control vector, likely in a local minimum, is a function of the starting control vector $[\theta_1^0,\theta_2^0,...\theta_N^0,\Omega^0]$. 

Instead, utilizing the ability of concurrent estimation to simulate many samples paths at once, we explore more of the control space by selecting control vectors randomly. Control vectors may be selected randomly using the discrete uniform distribution: 
\fbox{\begin{minipage}{23.75em}
\textbf{UNIFORMLY RANDOMLY SELECT A CONTROL PARAMETER VECTOR}

 1. SET CONTROL VECTOR 
 
 $\qquad \theta_i \sim$unif$\{0,m\} \forall i \in \mathcal{N}$ , $\Omega \sim$unif$\{1,m\} $

 2. DETERMINE IF CONTROL VECTOR FEASIBLE 
 
 $\qquad$ if $\sum_0^N \theta_i > m$ return to step 1. 

 \end{minipage}}
in which the second step verifies that the control vector is feasible. Note that if $\Omega > \sum_0^N \theta_i$, the control cannot be triggered rendering a feasible but clearly useless control. 

The control space of this AT load balancing problem is so large that it is not feasible to randomly simulate even a small percent of the possible controls for a long enough time to get a decent estimation of the performance. This justifies the following algorithm that finds families of well performing control vectors by iteratively shrinking the permissible control space. Instead of selecting control vectors by $\sim$unif$\{1,m\}$, let us define an lower and upper bound $\underline{\theta_i}$ and $\overline{\theta_i}$ or $\underline{\Omega}$ and $\overline{\Omega}$. These bounds will iteratively decrease the size of the search-able control space resulting in small permissible bands for each control parameter. Let the first iteration of the algorithm start with $\underline{\theta_i} = 0, \overline{\theta_i} = m \ \forall i \in \mathcal{N}, \underline{\Omega} = 1, \overline{\Omega} = m$. 

Let $T$ be the length of the simulated sample paths, $L$ be the number of sample paths to concurrently simulate, and $K$ be the number of times to run a $L$ concurrent sample paths:
\fbox{\begin{minipage}{23.75em}
\textbf{ITERATION OF RANDOM SEARCH ALGORITHM}

\textbf{for} $k$ = 1:$K$
 
$\quad$1. Find $L$ feasible control vectors.
 
$\quad$ $\quad$ \textbf{for} $l$ = 1:$L$
 
$\quad$ $\qquad$ 1A. Uniformly randomly select the $l$th control 

$\quad$ $\qquad$  vector: $\theta_i \sim$unif$\{\underline{\theta_i},\overline{\theta_i}\} \ \forall i \in \mathcal{N}$ , $\Omega \sim$unif$\{\underline{\Omega},\overline{\Omega}\} $  
 
$\quad$ $\qquad$ 1B. If $\sum_0^N \theta_i > m$, return to step 1A. 
 
$\quad$ $\quad$ \textbf{end}

$\quad$ 2. Simulate $L$ sample paths with control vectors from 

$\quad$ step 1 for $T$ time units via concurrent estimation.

$\quad$ 3. Choose the best performing control vector and label 

$\quad$ it $[\theta_1^{*k},...,\theta_N^{*k},\Omega^{*k}]$. 

\textbf{end}

\textbf{Redefine:} $\underline{\theta_i} = \min_k\{ \theta_i^{*k} \}$, $\overline{\theta_i} = \max_k\{ \theta_i^{*k} \}$,

$\underline{\Omega} = \min_k\{ \Omega^{*k} \}$, $\overline{\Omega} = \max_k\{ \Omega^{*k} \}$,

 \end{minipage}}

Each iteration of the above algorithm will decrease the size of the search-able control space as poorly performing control parameters are excluded in the \textbf{Redefine} step. This method, like the iterative search, may find areas of local minima. However, this algorithm may be run multiple times in independent trials to show that the final families of controls are very similar. 

The choice of $K,L$, and $T$ as well as the number of iterations will determine the final size and composition of the search-able control space. $T$, the length of the simulation, must be long enough to produce decent estimations of the control performance, but for the first couple of iterations it need not be too long as these first iterations weed out the especially poor performing control parameters. As the bands of search control parameters $[\underline{\theta_i},\overline{\theta_i}]$ thins, $T$ should increase such as to produce better estimations of system performance under the various control vectors. $L$, the number of sample paths concurrently estimated, should be large enough to render a good representation of the search space. Likewise, $K$ should be large enough that many possible control paths are simulated. 

Note that we could concurrently simulate $K \times L$ sample paths and pick the best $K$ control vectors to label as $[\theta_1^{*k},...,\theta_N^{*k},\Omega^{*k}], k = 1:K$. By instead choosing to separate them into $K$ series of $L$ sample paths, we mitigate the risk of an unusual sample path that disproportionately favors a select group of control vectors.

\section{System Performance Lower Bound} \label{Section_Lower_Bound}

In order to properly assess the performance of the parametric controllers
developed, we establish in what follows a lower bound on best performance
possible \textit{on average} (i.e., it is possible to construct a sample path
that performs better such as one without any $\delta_{i,j}$ request
events)with an objective function equal to the absolute lower bound of zero).
As the objective function in \eqref{J_intro} is an average over $[0,T]$, we
seek only an average performance lower bound.

Let us abstract the arrival process of discrete request events at rate
$\lambda_{i,j}$ into a continuous request flow process with rate
$\lambda_{i,j}$. On average, region $i$ will have an inflow of
$\sum_{j=1}^{N}\lambda_{j,i}$ and an outflow of $\sum_{j=1}^{N}\lambda_{i,j}$.
The average difference in request flow is given by: $d_{i}=\sum_{j=1}%
^{N}(\lambda_{i,j}-\lambda_{j,i}) $ such that we may define the following two
sets of regions depending on the sign of $d_{i}$:

$\mathcal{G}=\{i\in\mathcal{N}|d_{i}<0\} \qquad\mathcal{B}=\{i\in
\mathcal{N}|d_{i}\geq0\}$

Based on the objective function defined in \eqref{J_intro}, there are two
sources of system costs: unsated user requests for ATs and empty AT traveling
time. Regions within $\mathcal{B}$ will run out of AT flow as they are a
more popular origin than destination, thus they will effectively be forced to
reject a request flow $d_{i}$.
%As illustrated in Fig. \ref{Lost_Flow},
%consider a system with a request flow $\Lambda$ that is unable to sate some
%amount of flow $p$; thus the fraction of lost flow is $\frac{p}{\Lambda}$.
%\begin{figure}[th]
%\centering
%\includegraphics[width= 5 cm]{Lost_Flow.png}\caption{A system with inflow of
%$\Lambda$ and lost flow $p$ results in $\frac{p}{\Lambda}$ percent of flow
%lost.}%
%\label{Lost_Flow}%
%\end{figure}
%Likewise,

For this abstracted flow AT system an amount $p$ of flow lost by not being
sated with available ATs costs the system $c(p)$ as defined by:
\begin{equation}
c(p)=p\cdot\frac{1}{\sum_{i=1}^{N}\sum_{j=1}^{N}\lambda_{i,j}}
\label{cost_ignore}%
\end{equation}

Regions within $\mathcal{B}$ will build up excess AT flow and are
candidates to send empty AT flow out. Suppose that $p$ empty ATs are sent from
$i$ to $j$ at the beginning of a time period, with an average trip time
$\frac{1}{\mu_{i,j}}$ for a total mean empty AT driving time$\frac{p}%
{\mu_{i,j}}$. As there are a total of $m$ AT-hours, the mean number of empty
ATs driving is $\frac{p}{\mu_{i,j}m}$. Similarly, for the abstracted fluid AT
system forcing an AT flow of $p$ from $i$ to $j$, the system cost incurred is
as follows; Note that both \eqref{cost_ignore} and \eqref{cost_sate} are
linear functions of the flow $p$.
\begin{equation}
C(p,i,j)=p\cdot\frac{1}{\mu_{i,j}m} \label{cost_sate}%
\end{equation}

Consider two types of decision variables: $\beta_{j},\ j\in\mathcal{B}$, as
the fraction of positive request difference $d_{j}$ that will be left unsated
and $v_{i,j}$ as the forced empty AT flow from $i\in\mathcal{G}$ to
$j\in\mathcal{B}$ that will sate the remaining $[1-\beta_{j}]$ fractional
difference in flow. The following linear program finds a lower bound with less
than $N^{2}+N$ decision variables. The first and second parts of the objective
function are the cost of ignoring a fraction $\beta_{j}$ and sating the
fraction $[1-\beta_{j}]$ of difference in request $d_{j}$ from
\eqref{cost_ignore} and \eqref{cost_sate}, respectively. The first constraint
requires that the fraction $[1-\beta_{j}]$ of difference in request is sated
in \textquotedblleft bad\textquotedblright\ regions. The second
constraint places limitations on the empty AT flows from \textquotedblleft
good\textquotedblright\ regions.
\begin{equation}
\mbox{LB}=\min_{v_{i,j},\beta_{j}}\sum_{j\in\mathcal{B}}\bigg (\frac
{d_{j}\beta_{j}}{\sum_{i=1}^{N}\sum_{k=1}^{N}\lambda_{i,k}}+\sum
_{i\in\mathcal{G}}\frac{v_{i,j}}{\mu_{i,j}m}\bigg ) \label{LB_LP}%
\end{equation}%
\[
\begin{aligned}[t] \text{s. t.} \qquad d_j(1-\beta_j) &= \sum_{i \in \mathcal{G}} v_{i,j} & \quad j \in \mathcal{B} \\ & -d_i \ge \sum_{j \in \mathcal{B}} v_{i,j} & \quad i \in \mathcal{G} \\ & 0 \le v_{i,j} & \quad i \in \mathcal{G}, j \in \mathcal{B} \\ & 0 \le B_j \le 1 & \quad j \in \mathcal{B} \end{aligned}
\]

\section{Simulation Example: A 6-Region System}

\label{Example}

Consider an $N$=6 MoD system with request and travel rates as shown in Table
\ref{lambda} and an objective function weight $w$=0.5. In order to assess the
performance of our event-driven parametric controller, we compare it to the
time-driven controller, the static controller in \cite{pavone2012robotic}, the
lower bound derived in (\ref{LB_LP}), and the case of no control whatsoever.

%%
\begin{comment}
{\renewcommand{\arraystretch}{1.001}
\begin{table}[h]
\centering
\caption{N=6 System Request Rates $\lambda$}
\label{lambda}
\begin{tabular}{|l|l|l|l|l|l|}
\hline
6 & 15 & 6 & 6 & 9 & 3  \\ \hline
3 & 6  & 3 & 6 & 6 & 12 \\ \hline
0 & 9  & 3 & 0 & 3 & 3  \\ \hline
6 & 3  & 0 & 6 & 3 & 0  \\ \hline
6 & 12 & 6 & 0 & 3 & 0  \\ \hline
6 & 18 & 3 & 3 & 6 & 6  \\ \hline
\end{tabular}
\end{table}}
{\renewcommand{\arraystretch}{1.001}
\begin{table}[h]
\centering
\caption{N=6 System Travel Rates $\mu$}
\label{mu}
\begin{tabular}{|l|l|l|l|l|l|}
\hline
12  & 9.6 & 4.8 & 8.4 & 2.4 & 3.6 \\ \hline
9.6 & 12  & 7.2 & 6   & 3.6 & 4.8 \\ \hline
4.8 & 7.2 & 12  & 4.8 & 3.6 & 9.6 \\ \hline
8.4 & 6   & 4.8 & 12  & 2.4 & 2.4 \\ \hline
2.4 & 3.6 & 3.6 & 2.4 & 12  & 8.4 \\ \hline
3.6 & 4.8 & 9.6 & 2.4 & 8.4 & 12  \\ \hline
\end{tabular}
\end{table}}
%%
\end{comment}

\begin{table}[h]
\caption{$N$=6 System Request Rates $\lambda$ and Travel Rates $\mu$}%
\label{lambda}
\centering
\begin{tabular}
[c]{ccccccccccccc}%
\multicolumn{6}{c}{$\lambda$ (demands/min)} &  & \multicolumn{6}{c}{$\mu$
(1/min)}\\\cline{1-6}\cline{8-13}%
\multicolumn{1}{|c|}{$\!$6$\!$} & \multicolumn{1}{c|}{$\!$15$\!$} &
\multicolumn{1}{c|}{6} & \multicolumn{1}{c|}{6} & \multicolumn{1}{c|}{9} &
\multicolumn{1}{c|}{3} & \multicolumn{1}{c|}{} & \multicolumn{1}{c|}{12} &
\multicolumn{1}{c|}{9.6} & \multicolumn{1}{c|}{4.8} & \multicolumn{1}{c|}{8.4}
& \multicolumn{1}{c|}{2.4} & \multicolumn{1}{c|}{3.6}\\\cline{1-6}\cline{8-13}%
\multicolumn{1}{|c|}{$\!$3$\!$} & \multicolumn{1}{c|}{$\!$6$\!$} &
\multicolumn{1}{c|}{3} & \multicolumn{1}{c|}{6} & \multicolumn{1}{c|}{6} &
\multicolumn{1}{c|}{$\!$12$\!$} & \multicolumn{1}{c|}{} &
\multicolumn{1}{c|}{9.6} & \multicolumn{1}{c|}{12} & \multicolumn{1}{c|}{7.2}
& \multicolumn{1}{c|}{6} & \multicolumn{1}{c|}{3.6} & \multicolumn{1}{c|}{4.8}%
\\\cline{1-6}\cline{8-13}%
\multicolumn{1}{|c|}{$\!$0$\!$} & \multicolumn{1}{c|}{9} &
\multicolumn{1}{c|}{3} & \multicolumn{1}{c|}{0} & \multicolumn{1}{c|}{3} &
\multicolumn{1}{c|}{3} & \multicolumn{1}{c|}{} & \multicolumn{1}{c|}{4.8} &
\multicolumn{1}{c|}{7.2} & \multicolumn{1}{c|}{12} & \multicolumn{1}{c|}{4.8}
& \multicolumn{1}{c|}{3.6} & \multicolumn{1}{c|}{9.6}\\\cline{1-6}\cline{8-13}%
\multicolumn{1}{|c|}{$\!$6$\!$} & \multicolumn{1}{c|}{3} &
\multicolumn{1}{c|}{0} & \multicolumn{1}{c|}{6} & \multicolumn{1}{c|}{3} &
\multicolumn{1}{c|}{0} & \multicolumn{1}{c|}{} & \multicolumn{1}{c|}{8.4} &
\multicolumn{1}{c|}{6} & \multicolumn{1}{c|}{4.8} & \multicolumn{1}{c|}{12} &
\multicolumn{1}{c|}{2.4} & \multicolumn{1}{c|}{2.4}\\\cline{1-6}\cline{8-13}%
\multicolumn{1}{|c|}{$\!$6$\!$} & \multicolumn{1}{c|}{$\!$12$\!$} &
\multicolumn{1}{c|}{6} & \multicolumn{1}{c|}{0} & \multicolumn{1}{c|}{3} &
\multicolumn{1}{c|}{0} & \multicolumn{1}{c|}{} & \multicolumn{1}{c|}{2.4} &
\multicolumn{1}{c|}{3.6} & \multicolumn{1}{c|}{3.6} & \multicolumn{1}{c|}{2.4}
& \multicolumn{1}{c|}{12} & \multicolumn{1}{c|}{8.4}\\\cline{1-6}\cline{8-13}%
\multicolumn{1}{|c|}{$\!$6$\!$} & \multicolumn{1}{c|}{$\!$18$\!$} &
\multicolumn{1}{c|}{3} & \multicolumn{1}{c|}{3} & \multicolumn{1}{c|}{6} &
\multicolumn{1}{c|}{6} & \multicolumn{1}{c|}{} & \multicolumn{1}{c|}{3.6} &
\multicolumn{1}{c|}{4.8} & \multicolumn{1}{c|}{9.6} & \multicolumn{1}{c|}{2.4}
& \multicolumn{1}{c|}{8.4} & \multicolumn{1}{c|}{12}\\\cline{1-6}\cline{8-13}%
\end{tabular}
\end{table}

The event-driven and time-driven controllers require $N+1$ integer and $1$
real-valued parameters respectively, which are determined using the concurrent
estimation techniques
described in section \ref{SC}. After running many iterations of the greedy iteration as described by Section \ref{Greedy_Iterative_Search_AT} on a shared cloud computer cluster in
MATLAB 2018b, the event and time-driven parameters in Table
\ref{Driven_Parameters} were determined to perform the best for fleet sizes of
50, 75, 100, and 125.
%over a simulated time of $T=100,000$ minutes
\makeatletter
\def\hlinewd#1{%
  \noalign{\ifnum0=`}\fi\hrule \@height #1 \futurelet
   \reserved@a\@xhline}
\makeatother
\newcolumntype{?}{!{\vrule width 1.5pt}}
{\renewcommand{\arraystretch}{1.00}
\begin{table}[h]
\centering
\caption{Event and Time-Driven Controller Parameters: $N$=6 System}
\begin{tabular}{|c?c|c|c|c|c|c|c?c|}
\hline
\textbf{Control} & \multicolumn{7}{c?}{\textbf{Event-Driven}} & \textbf{Time-Driven} \\ \hline
$m$ & $\theta_1$ & $\theta_2$ & $\theta_3$ & $\theta_4$ & $\theta_5$ & $\theta_6$ & $\Omega$ & $\Omega$  \\  \hlinewd{1.5pt} 
50 & 10	& 7	& 4	& 1	& 4	& 7	& 5 & 24\\ \hline
75 & 15	& 13 & 8 & 4 & 12 & 13 & 8 & 12 \\ \hline
100 & 20& 	16	& 11& 	7& 	16& 	20& 	14 & 12\\ \hline
125 & 27& 	19& 	13& 	7& 	19& 	25& 	22 & 18\\ \hline
\end{tabular}
\label{Driven_Parameters}
\end{table}}

As these parameter vectors were determined to be the best via simulation, they
are likely local minima of the objective function. Iterations stopped when
slight deviations (i.e. $\theta_{i}^{\prime}=\theta_{i}+1$) had little effect
on the objective function of a sufficiently long sample path ($T$=100,000 time
units). 

Likewise, families of control vectors were found using the alternate random search method of finding well performing vectors as described in Section \ref{Random_Search_AT} for each fleet size of [50,75,100,125] with iteration-specific $T,L$, and $K$: 
\begin{table}[]
\centering
\caption{Random Search Iteration-Specific Parameters}
\label{tab:my-table}
\begin{tabular}{|c|c|c|c|c|c|}
\hline
Iteration & 1-3 & 4-6   & 7-9    & 10-12    & 13-15   \\ \hline
$T$       & 500 & 5,000 & \multicolumn{3}{c|}{10,000} \\ \hline
$L$       & 25  & 50    & 100    & 200      & 500     \\ \hline
$K$       & \multicolumn{5}{c|}{25}                   \\ \hline
\end{tabular}
\end{table}

Table \ref{Random_Search_Table} shows the final well performing families of control vectors from two independent trials of 15 successive iterations. Note that these independent trials have very similar final control vector families and include the control vectors found via iterative greedy search as in Table \ref{Driven_Parameters}. 
\begin{table*}[!ht]
\centering
\caption{Families of Well Performing Vectors: Iterative Random Search}
\label{Random_Search_Table}
\begin{tabular}{|c|c|c|c|c|c|c|c|}
\hline
$m$                  & $[\underline{\theta_1},\overline{\theta_1}]$ & $[\underline{\theta_2},\overline{\theta_2}]$ & $[\underline{\theta_3},\overline{\theta_3}]$ & $[\underline{\theta_4},\overline{\theta_4}]$ & $[\underline{\theta_5},\overline{\theta_5}]$ & $[\underline{\theta_6},\overline{\theta_6}]$ & $[\underline{\Omega},\overline{\Omega}]$ \\ \hline
\multirow{2}{*}{50}  & {[}9,12{]}                                   & {[}4,9{]}                                    & {[}2,5{]}                                    & {[}0,3{]}                                    & {[}1,5{]}                                    & {[}6,8{]}                                    & {[}4,8{]}                                \\ \cline{2-8} 
                     & {[}9,12{]}                                   & {[}1,9{]}                                    & {[}2,5{]}                                    & {[}0,3{]}                                    & {[}0,5{]}                                    & {[}6,8{]}                                    & {[}3,8{]}                                \\ \hline
\multirow{2}{*}{75}  & {[}15,16{]}                                  & {[}12,14{]}                                  & {[}7,8{]}                                    & {[}2,5{]}                                    & {[}10,13{]}                                  & {[}13,15{]}                                  & {[}8,11{]}                               \\ \cline{2-8} 
                     & {[}14,17{]}                                  & {[}11,14{]}                                  & {[}7,9{]}                                    & {[}0,6{]}                                    & {[}10,13{]}                                  & {[}12,15{]}                                  & {[}6,11{]}                               \\ \hline
\multirow{2}{*}{100} & {[}19,23{]}                                  & {[}15,17{]}                                  & {[}9,12{]}                                   & {[}4,8{]}                                    & {[}14,18{]}                                  & {[}18,21{]}                                  & {[}12,17{]}                              \\ \cline{2-8} 
                     & {[}18,22{]}                                  & {[}15,18{]}                                  & {[}9,12{]}                                   & {[}1,8{]}                                    & {[}14,17{]}                                  & {[}17,22{]}                                  & {[}10,16{]}                              \\ \hline
\multirow{2}{*}{125} & {[}24,30{]}                                  & {[}17,21{]}                                  & {[}11,15{]}                                  & {[}5,12{]}                                   & {[}18,22{]}                                  & {[}22,28{]}                                  & {[}16,25{]}                              \\ \cline{2-8} 
                     & {[}24,30{]}                                  & {[}16,22{]}                                  & {[}11,16{]}                                  & {[}2,11{]}                                   & {[}17,21{]}                                  & {[}23,28{]}                                  & {[}16,26{]}                              \\ \hline
\end{tabular}
\end{table*}

The static controller introduced in \cite{pavone2012robotic} requires
$N^{2}$ parameters found via the LP \eqref{MINLP} and shown in Table
\ref{Static_Parameters}. 
{\renewcommand{\arraystretch}{1}
\begin{table}[h]
\caption{Static Controller Parameters for the $N$=6 System}
\centering
\begin{tabular}{cc?c|c|c|c|c|c|}
\cline{3-8}
 &  & \multicolumn{6}{c|}{Destination} \\ \cline{3-8} 
 & $r_{i,j}$ & 1 & 2 & 3 & 4 & 5 & 6 \\ \hlinewd{1.5pt} 
\multicolumn{1}{|c|}{\multirow{6}{*}{Origin}} & 1 & 0 & 0 & 0 & 0 & 0 & 0 \\ \cline{2-8} 
\multicolumn{1}{|c|}{}  & 2 & 0.25 & 0 & 0 & 0 & 0 & 0.2 \\ \cline{2-8} 
\multicolumn{1}{|c|}{} & 3 & 0 & 0 & 0 & 0 & 0 & 0.05 \\ \cline{2-8} 
\multicolumn{1}{|c|}{} & 4 & 0.5 & 0 & 0 & 0 & 0 & 0 \\ \cline{2-8} 
\multicolumn{1}{|c|}{} & 5 & 0 & 0 & 0 & 0 & 0 & 0.05 \\ \cline{2-8} 
\multicolumn{1}{|c|}{} & 6 & 0 & 0 & 0 & 0 & 0 & 0 \\ \hline
\end{tabular}
\label{Static_Parameters}
\end{table}}

\begin{figure}[h]
\centering
\includegraphics[width=1\columnwidth]{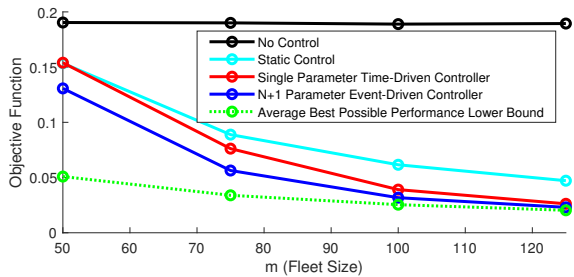}\caption{The
performance of a system with fleet sizes 50, 75, 100, and 125 under no
control, static control, time-driven control, event-driven control, and the
average best possible performance lower bound.}%
\label{compare_obj}%
\end{figure}

Figure \ref{compare_obj} shows the average simulated performance for all fleet
sizes and controllers.
%for a simulation time $T=1,000,000$ minutes \textbf{[XXX See comment on time units]}.
All systems performed about the same under no control -- over 37$\%$ of user
requests unsated. As expected, the $N$+1 parameter event-driven controller
with its state-dependent control and system-specific tuned parameters performs
the best across all fleet sizes.

The intermediate fleet sizes studied here are where the true benefits of the
event-driven controller are to be observed for this particular 6-region
system. Note that the event and time-driven controllers quickly approach the
lower bound as the fleet size increases but the static controller does not--
this advantage is due to the event and time-driven controllers use of state
information. We ignore fleets of less than 50: they are unstable and perform
poorly no matter the control, since the underlying MoD system is
under-capacitated with an insufficient number of vehicles to satisfy the given
demand. We also ignore fleets over 125: they perform well no matter the
control, as they are over-capacitated.

Table \ref{Objective_75} shows a detailed performance comparison of the three
controllers relative to the uncontrolled case and the lower bound derived in
(\ref{LB_LP}) for a fleet size of 75 ATs.
{\renewcommand{\arraystretch}{1.0}
\begin{table}[h]
\centering
\caption{$N$=6, $m$=75 System Performance Comparison}
\begin{tabular}{|c|c|c|c|c|c|}
\hline
Control & None & Static & \begin{tabular}[c]{@{}c@{}}Time\\ Driven\end{tabular} & \begin{tabular}[c]{@{}c@{}}Event\\ Driven\end{tabular} & \begin{tabular}[c]{@{}c@{}}Lower \\ Bound\end{tabular} \\ \hline
\% Users Rejected & 38.0  & 11.8 & 7.0 & 3.4 & 0 \\ \hline
\% ATs Drive Empty & 0  & 6.0 & 8.2 & 7.8 & 6.8 \\ \hline
 $\bar{J} \ \ (w=0.5)$ & 19.0	& 8.9 &	7.6 &	5.6 &	3.4 \\ \hline
\end{tabular}
\label{Objective_75}
\end{table}}

Table \ref{Infeasible_Events_Example}, a complement to Table \ref{Infeasible_Events}, shows the range of percents of fictitious events created by the Standard Clock method and the variation thereof in Section \ref{SC} for a fleet size of $m=75$ ATs. The variation of the Standard Clock method avoids many more fictitious events -- enough to justify the extra calculations required. 

{\renewcommand{\arraystretch}{1.1}
\begin{table}[h]
\centering
\caption{Percent Fictitious Events for Fleet Size $m=75$}
\label{Ficticious_Events_Example}
\begin{tabular}{c|c|c|}
\cline{2-3}
                                     & Minimum & Maximum \\ \hline
\multicolumn{1}{|l|}{Standard Clock} &  93.91  &   98.96  \\ \hline
\multicolumn{1}{|l|}{Variation}      &   0    &  82.87   \\ \hline
\end{tabular}
\label{Infeasible_Events_Example}
\end{table}}

\section{Conclusions and Future Work}

\label{Conclusions}

MoD systems must load balance by sending empty vehicles to mitigate the
temporal demand patterns that deplete some service regions of available
vehicles. We have defined an objective function to jointly minimize the
fraction of user requests denied due to unavailability and the fraction of
time vehicles drive empty-- and derived its lower bound. As optimal control
via DP quickly becomes intractable even for small dimensionality systems, we
have developed a parametric controller using thresholds on the number of
vehicles available in and en route to each region. Optimal (or at least
well-performing) parameters were determined using Concurrent Estimation
methods which allow for the construction of multiple sample paths under
different parameters from a single nominal sample path. Our simulation
examples for the proposed event-driven threshold-based controller perform
significantly better than an uncontrolled system or static controller and
approach the lower bound for large fleet sizes.
%In order to find well performing control policy parameters we introduce a variation on the Standard Clock method of concurrent estimation which takes advantage of MoD system framework to hasten sample path construction.

Future work will include using actual taxi data and exploring a wider range of
the control parameter space by making a more efficient use of CE methods
through which we may attain a global optimum.

\bibliographystyle{IEEEtran}
\bibliography{library}

\section*{Appendix A: Optimal Control Policies of Simple Systems } \label{A1}

\begin{figure}[th]
\centering
\includegraphics[width= 7 cm]{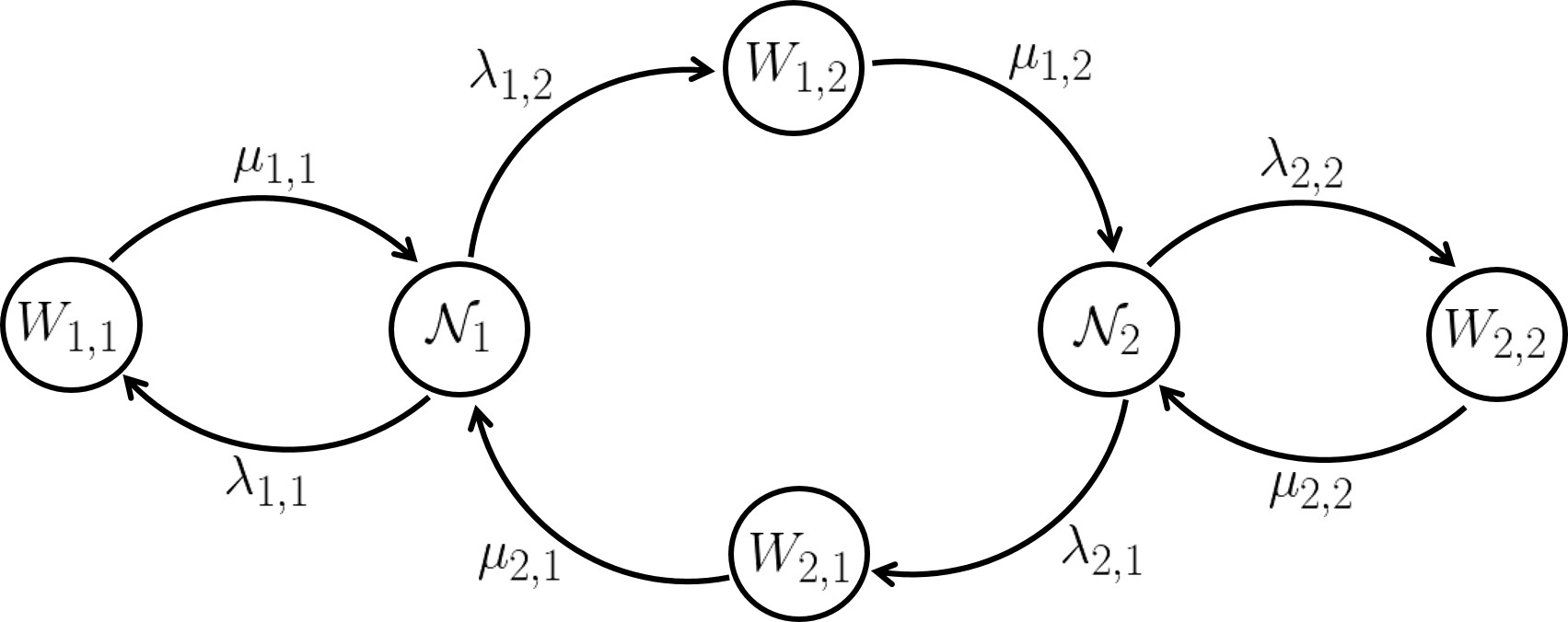} \caption{The simplest possible system, the time-invariant 2-region 1-AT system sans controls has six possible states: the AT may be in either of the idle AT queues $\mathcal{N}_i$ or in one of the four servers $W_{i,j}$; the system moves between these six states at rates $\mu_{i,j}$ and $\lambda_{i,j}$  where $\lambda_{i,j} = p_{i,j}\lambda_{i}$.  }%
\label{chain}
\end{figure}

\begin{figure}[th]
\centering
\includegraphics[width= 7 cm]{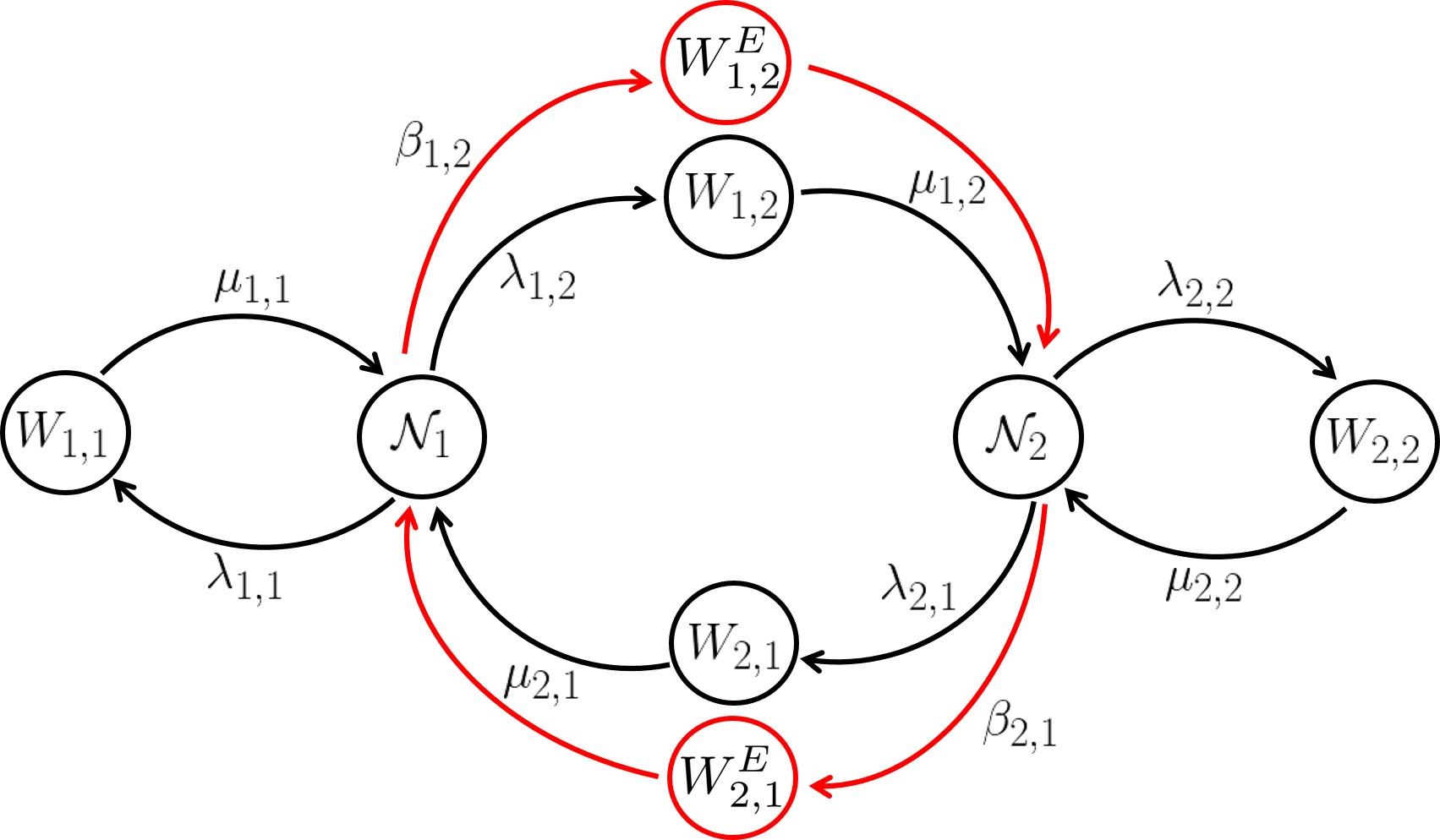} \caption{The time-invariant $N=2$, $m=1$ system with controllable transition rates $\beta_{1,2}$ and $\beta_{2,1}$ as the rates at which empty ATs are sent from 1 to 2 and 2 to 1, respectively, leads to an 8 state system with the addition of empty AT states $W_{1,2}^E$ and $W_{2,1}^E$.}%
\label{chain_control}
\end{figure}

The optimal control policies for simple systems are possible to determine because they have tractable steady-state distributions.
For a simple 2-region 1-AT time-invariant system (omit $k$) there are only eight possible states: the single AT is either in one of the two idle AT queues $\mathcal{N}_i$, full with a passenger in server $W_{i,j}$ (including $i = j$), or empty in server $W_{i,j}$ (excluding the obviously undesirable empty intra-region trips $i = j$) as listed in Table \ref{simple_states}.
\begin{table}[h]
\centering
\caption{States of the Time-Invariant $N=2,m=1$ System}
\label{simple_states}
\begin{tabular}{|c|*3{>{\renewcommand{\arraystretch}{1}}c}|}
\hline
State & $\left[ \begin{array}{cc} x_1 & x_2  \end{array}\right]$ & $\left[ \begin{array}{cc} y_{1,1} & y_{1,2}  \\ y_{2,1} & y_{2,2} \end{array}\right]$ & $\left[ \begin{array}{cc} z_{1,1} & z_{1,2} \\ z_{2,1} & z_{2,2}  \end{array}\right]$\\
\hline
$N_1$ & $\left[ \begin{array}{cc} 1 & 0  \end{array}\right]$ & $\left[ \begin{array}{cc} 0 & 0  \\ 0 & 0 \end{array}\right]$ & $\left[ \begin{array}{cc} 0 & 0 \\ 0 & 0  \end{array}\right]$\\
\hline
$N_2$ & $\left[ \begin{array}{cc} 0 & 1  \end{array}\right]$ & $\left[ \begin{array}{cc} 0 & 0  \\ 0 & 0 \end{array}\right]$ & $\left[ \begin{array}{cc} 0 & 0 \\ 0 & 0  \end{array}\right]$\\
\hline
$W_{1,1}$ & $\left[ \begin{array}{cc} 0 & 0  \end{array}\right]$ & $\left[ \begin{array}{cc} 1 & 0  \\ 0 & 0 \end{array}\right]$ & $\left[ \begin{array}{cc} 0 & 0 \\ 0 & 0  \end{array}\right]$\\
\hline
$W_{1,2}$ & $\left[ \begin{array}{cc} 0 & 0  \end{array}\right]$ & $\left[ \begin{array}{cc} 0 & 1  \\ 0 & 0 \end{array}\right]$ & $\left[ \begin{array}{cc} 0 & 0 \\ 0 & 0  \end{array}\right]$\\
\hline
$W_{2,1}$ & $\left[ \begin{array}{cc} 0 & 0  \end{array}\right]$ & $\left[ \begin{array}{cc} 0 & 0  \\ 1 & 0 \end{array}\right]$ & $\left[ \begin{array}{cc} 0 & 0 \\ 0 & 0  \end{array}\right]$\\
\hline
$W_{2,2}$ & $\left[ \begin{array}{cc} 0 & 0  \end{array}\right]$ & $\left[ \begin{array}{cc} 0 & 0  \\ 0 & 1 \end{array}\right]$ & $\left[ \begin{array}{cc} 0 & 0 \\ 0 & 0  \end{array}\right]$\\
\hline
$W_{1,2}^E$ & $\left[ \begin{array}{cc} 0 & 0  \end{array}\right]$ & $\left[ \begin{array}{cc} 0 & 0  \\ 0 & 0 \end{array}\right]$ & $\left[ \begin{array}{cc} 0 & 1 \\ 0 & 0  \end{array}\right]$\\
\hline
$W_{2,1}^E$ & $\left[ \begin{array}{cc} 0 & 0  \end{array}\right]$ & $\left[ \begin{array}{cc} 0 & 0  \\ 0 & 0 \end{array}\right]$ & $\left[ \begin{array}{cc} 0 & 0 \\ 1 & 0  \end{array}\right]$\\
\hline
\end{tabular}
\end{table}
Let the user arrivals follow a Poisson process with rate $\lambda_{i}$ and infinite-capacity server $W_{i,j}$ be exponential with mean service time $\frac{1}{\mu_{i,j}}$.
The simple $N=2,m=1$ system without controls (omit empty driving states $W_{i,j}^E$) may now be represented by the continuous Markov chain in Figure \ref{chain}. For simplicity of notation, let $\lambda_{i,j} = p_{i,j}\lambda_{i}$, the rate at which users arrive to region $i$ destined for region $j$.
Instead of individual AT controls $u_{i,j,q}$, consider a controllable state transition rate $\beta_{i,j} \in [0,\infty)$ as the \textit{rate} at which to send an empty AT from region $i$ to region $j$ \cite{zhang2016control},\cite{spieser2014toward}. This controllable Markov Chain for the $N=2,m=1$ system in Figure \ref{chain_control} includes empty AT states $W_{1,2}^E$ and $W_{2,1}^E$.
We may calculate the steady-state probability vector $\pmb{\pi}_{\pmb{\beta}}$ from the normalized transition matrix with rates from Figure \ref{chain_control} where $\pmb{\beta} = [\beta_{1,2},\beta_{2,1}]$.
The cost $C(\mathscr{X})$ associated with being in state $\mathscr{X}$ is as follows:
\begin{equation}
C(\mathscr{X}) = w \sum_{i=1}^N \frac{ \sum_j \lambda_{i,j} \mathbf{1}[x_i = 0] }{ \sum_i \sum_j (\lambda_{i,j} +\mu_{i,j}) } + (1-w) \frac{ \sum_i \sum_j z_{i,j} } { m }
\label{state_x_cost_Appendix}
\end{equation}
The full objective function of this $N=2,m=1$ system may be written in terms of the probabilities:
\begin{multline}
J^* \! \! = \! \! \min_{\pmb{\beta}}  \!  \bigg \{ \! \! w  \! \frac{\left(\lambda_{1,1} \! + \! \lambda_{1,2}\right) \! \! (1 \text{-} \pi_{\pmb{\beta}} \! (\mathcal{N}_1 \! )  \! ) \! + \! \left(\lambda_{2,1} \! + \! \lambda_{2,2}\right) \! \! (1 \text{-} \pi_{\pmb{\beta}} \! (\mathcal{N}_2 \! ) \! )}{\sum_{i}\sum_j \lambda_{i,j}}   \\
+ (1-w)  \big ( \pi_{\pmb{\beta}}(W^E_{1,2})  + \pi_{\pmb{\beta}}(W^E_{2,1})  \big ) \bigg \}
\label{simple_obj}
\end{multline}
Consider the simplest scenario in which the percents of rejected requests and load balancing vehicles are equally weighted with $w=0.5$ such that we may ignore $w$. The objective function may be alternatively shown to be a ratio of polynomial functions of the $\beta_{1,2}$ and $\beta_{2,1}$ control rates:
\begin{equation}
f(\beta_{1,2},\beta_{2,1}) = \frac{A\beta_{1,2}\beta_{2,1} + B\beta_{1,2} + C\beta_{2,1} + D}{E\beta_{1,2}\beta_{2,1} + F\beta_{1,2} + G\beta_{2,1} + H}
\label{ABCs}
\end{equation}
where $A$ through $H$ are all positive complicated combinations of $\lambda_{i,j}$ and $\mu_{i,j}$ listed in Appendix B. 

Inspection of A and E in Appendix B reveals that $A=2E$ such that the objective function can be rewritten as:

\begin{equation}
f(\beta_{1,2},\beta_{2,1}) = \frac{A\beta_{1,2}\beta_{2,1} + B\beta_{1,2} + C\beta_{2,1} + D}{0.5A\beta_{1,2}\beta_{2,1} + F\beta_{1,2} + G\beta_{2,1} + H}
\label{ABCs2}
\end{equation}

As both $\beta_{1,2}$ and $\beta_{2,1} \rightarrow \infty$ the objective function goes to its upper bound of 2 as the AT cycles empty between the two regions and never serves any users:
\begin{equation}
\mbox{lim}_{\beta_{1,2}.\beta_{2,1} \rightarrow \infty} \ \ f(\beta_{1,2},\beta_{2,1}) = \frac{A}{.5A}=2
\end{equation}
leading to 100$\%$ of user requests missed and 100$\%$ of AT timing driving empty. 

The derivative of the objective function with regard to $\beta_{1,2}$
\begin{multline}
 \frac{\partial f(\beta_{1,2},\beta_{2,1}) }{\partial \beta_{1,2}} = \\
\frac{A\beta_{2,1} \! ( \! G\beta_{2,1} \! \! + \! H) \! + \! B( \! G\beta_{2,1} \! + \! H) \! + \! ( \! C\beta_{2,1} \! + \! D \! )( \! \text{-} \! E\beta_{2,1} \text{-} F)}{{\left(H+F\,\beta _{1,2}+G\,\beta _{2,1}+E\,\beta _{1,2}\,\beta _{2,1}\right)}^2}
\end{multline}
has an either always positive or negative numerator which is \emph{not} a function of $\beta_{1,2}$ and a denominator that is an always positive function of $\beta_{1,2}$. Hence $ \frac{\partial f(\beta_{1,2},\beta_{2,1}) }{\partial \beta_{1,2}}$ is either always positive or always negative such that the $\beta_{1,2}$ to minimize the objective function will be found on the boundaries 0 or $\infty$. 

Fixing $\beta_{2,1}$ to consider the objective as a function of $\beta_{1,2}$:
\begin{equation}
f(\beta_{1,2}) = \frac{(A\beta_{2,1} + B)\beta_{1,2} + (C\beta_{2,1} + D)}{(E\beta_{2,1} + F)\beta_{1,2} + (G\beta_{2,1} + H)} 
\end{equation}
we may regroup: $a = A\beta_{2,1} + B, \ b = C\beta_{2,1} + D, \ c = E\beta_{2,1} + F, \ d = G\beta_{2,1} + H$ with $a,b,c,d \ge 0$ to  rewrite the objective function as:
\begin{equation}
f(\beta_{1,2}) = \frac{a\beta_{1,2} + b}{c\beta_{1,2} +d} 
\end{equation}
where the first and second derivatives are:
\begin{equation}
\frac{d f(\beta_{1,2}) }{d \beta_{1,2}} = \frac{ad - cb}{(c\beta_{1,2}+d)^2}
\end{equation}
\begin{equation}
\frac{d^2 f(\beta_{1,2}) }{d \beta_{1,2}\ ^2} = \frac{-2c(c\beta_{1,2}+d)(ad - cb)}{(c\beta_{1,2}+d)^4}
\end{equation}

If we assume $ad > cb$ then the first derivative is always positive meaning the objective function is monotonically increasing and therefore minimized at $\beta_{1,2} = 0$ for $f(0) = \frac{b}{d}$. Furthermore the second derivative will always be negative implying the objective function is concave. Likewise if we assume $ad < cb$ then the first derivative is always negative meaning the objective function monotonically decreasing and is minimized at $\beta_{1,2} = \infty$ for $f(\infty) = \frac{a}{c}$. The second derivative is always positive and implies the objective function is convex. In the same way by fixing $\beta_{1,2}$, we may show that the optimal choice of $\beta_{2,1}$ is either at $0$ or $\infty$. 

% Alternately, we can show that $\beta_{1,2}$ and $\beta_{2,1}$ are either 0 or $\infty$ by the following logic.  ADD PROOF HERE

As the optimal choices for $\beta_{1,2}$ and $\beta_{2,1}$ are both either $0$ or $\infty$ there are only four possible sets of optimal controls and a simple comparison rule to determine the optimal policy:
\begin{equation}
(\beta_{1,2} , \beta_{2,1} ) = %
\begin{cases}
( \infty , \infty ) & \mbox{if} \min \{ \frac{A}{E} , \frac{B}{F} , \frac{C}{G} , \frac{D}{H} \} = \frac{A}{E} \\
( \infty , 0 ) & \mbox{if} \min \{ \frac{A}{E} , \frac{B}{F} , \frac{C}{G} , \frac{D}{H} \} = \frac{B}{F} \\
( 0 , \infty ) & \mbox{if} \min \{ \frac{A}{E} , \frac{B}{F} , \frac{C}{G} , \frac{D}{H} \} = \frac{C}{G} \\
( 0 , 0 ) & \mbox{if} \min \{ \frac{A}{E} , \frac{B}{F} , \frac{C}{G} , \frac{D}{H} \} = \frac{D}{H} \\
\end{cases}
\end{equation}
As previously noted, the first case will never occur as this would send the objective function to its upper bound of 2.

\section*{Appendix B: Optimal Policy for a $N=2,m=1$ Equally Weighted $w = 0.5$ System Using Rate $\beta_{i,j}$}\label{A2}

\[
A = 2\left(\mu _{11}\,\mu _{12}\,\mu _{22}+\mu _{11}\,\mu _{21}\,\mu _{22}\right)\,\left(\lambda _{11}+\lambda _{12}+\lambda _{21}+\lambda _{22}\right)
\]
\[
B =
\begin{array}{l}
{\lambda _{21}}^2\,\mu _{11}\,\mu _{12}\,\mu _{22}+{\lambda _{22}}^2\,\mu _{11}\,\mu _{12}\,\mu _{21}\\
+2\,{\lambda _{21}}^2\,\mu _{11}\,\mu _{21}\,\mu _{22}+\lambda _{11}\,\lambda _{21}\,\mu _{11}\,\mu _{12}\,\mu _{22} \\
+\lambda _{11}\,\lambda _{22}\,\mu _{11}\,\mu _{12}\,\mu _{21}
+\lambda _{12}\,\lambda _{21}\,\mu _{11}\,\mu _{12}\,\mu _{22}\\
+\lambda _{12}\,\lambda _{22}\,\mu _{11}\,\mu _{12}\,\mu _{21}
+2\,\lambda _{11}\,\lambda _{21}\,\mu _{11}\,\mu _{21}\,\mu _{22} \\
+2\,\lambda _{12}\,\lambda _{21}\,\mu _{11}\,\mu _{21}\,\mu _{22} 
+\lambda _{21}\,\lambda _{22}\,\mu _{11}\,\mu _{12}\,\mu _{21}\\
+\lambda _{21}\,\lambda _{22}\,\mu _{11}\,\mu _{12}\,\mu _{22}
+2\,\lambda _{21}\,\lambda _{22}\,\mu _{11}\,\mu _{21}\,\mu _{22} \\
+\lambda _{11}\,\mu _{11}\,\mu _{12}\,\mu _{21}\,\mu _{22}
+\lambda _{12}\,\mu _{11}\,\mu _{12}\,\mu _{21}\,\mu _{22}\\
\end{array}
\]
\[
C = 
\begin{array}{l}
2\,{\lambda _{12}}^2\,\mu _{11}\,\mu _{12}\,\mu _{22}
+{\lambda _{11}}^2\,\mu _{12}\,\mu _{21}\,\mu _{22}\\
+{\lambda _{12}}^2\,\mu _{11}\,\mu _{21}\,\mu _{22}
+2\,\lambda _{11}\,\lambda _{12}\,\mu _{11}\,\mu _{12}\,\mu _{22} \\
+\lambda _{11}\,\lambda _{12}\,\mu _{11}\,\mu _{21}\,\mu _{22}
+\lambda _{11}\,\lambda _{12}\,\mu _{12}\,\mu _{21}\,\mu _{22}\\
+2\,\lambda _{12}\,\lambda _{21}\,\mu _{11}\,\mu _{12}\,\mu _{22}
+2\,\lambda _{12}\,\lambda _{22}\,\mu _{11}\,\mu _{12}\,\mu _{22} \\
+\lambda _{11}\,\lambda _{21}\,\mu _{12}\,\mu _{21}\,\mu _{22}
+\lambda _{12}\,\lambda _{21}\,\mu _{11}\,\mu _{21}\,\mu _{22}\\
+\lambda _{11}\,\lambda _{22}\,\mu _{12}\,\mu _{21}\,\mu _{22}
+\lambda _{12}\,\lambda _{22}\,\mu _{11}\,\mu _{21}\,\mu _{22} \\
+\lambda _{21}\,\mu _{11}\,\mu _{12}\,\mu _{21}\,\mu _{22}
+\lambda _{22}\,\mu _{11}\,\mu _{12}\,\mu _{21}\,\mu _{22}
\end{array}
\]
\[
D = 
\begin{array}{l}
\lambda _{12}\,{\lambda _{21}}^2\,\mu _{11}\,\mu _{12}\,\mu _{22}
+\lambda _{12}\,{\lambda _{22}}^2\,\mu _{11}\,\mu _{12}\,\mu _{21}\\
+{\lambda _{12}}^2\,\lambda _{21}\,\mu _{11}\,\mu _{12}\,\mu _{22}
+{\lambda _{12}}^2\,\lambda _{22}\,\mu _{11}\,\mu _{12}\,\mu _{21} \\
+\lambda _{11}\,{\lambda _{21}}^2\,\mu _{12}\,\mu _{21}\,\mu _{22}
+\lambda _{12}\,{\lambda _{21}}^2\,\mu _{11}\,\mu _{21}\,\mu _{22}\\
+{\lambda _{11}}^2\,\lambda _{21}\,\mu _{12}\,\mu _{21}\,\mu _{22}
+{\lambda _{12}}^2\,\lambda _{21}\,\mu _{11}\,\mu _{21}\,\mu _{22} \\
+{\lambda _{12}}^2\,\mu _{11}\,\mu _{12}\,\mu _{21}\,\mu _{22}
+{\lambda _{21}}^2\,\mu _{11}\,\mu _{12}\,\mu _{21}\,\mu _{22}\\
+\lambda _{11}\,\lambda _{12}\,\lambda _{21}\,\mu _{11}\,\mu _{12}\,\mu _{22}
+\lambda _{11}\,\lambda _{12}\,\lambda _{22}\,\mu _{11}\,\mu _{12}\,\mu _{21} \\
+\lambda _{11}\,\lambda _{12}\,\lambda _{21}\,\mu _{11}\,\mu _{21}\,\mu _{22}+\lambda _{11}\,\lambda _{12}\,\lambda _{21}\,\mu _{12}\,\mu _{21}\,\mu _{22}\\
+\lambda _{12}\,\lambda _{21}\,\lambda _{22}\,\mu _{11}\,\mu _{12}\,\mu _{21}
+\lambda _{12}\,\lambda _{21}\,\lambda _{22}\,\mu _{11}\,\mu _{12}\,\mu _{22} \\
+\lambda _{11}\,\lambda _{21}\,\lambda _{22}\,\mu _{12}\,\mu _{21}\,\mu _{22}+\lambda _{12}\,\lambda _{21}\,\lambda _{22}\,\mu _{11}\,\mu _{21}\,\mu _{22}\\
+\lambda _{11}\,\lambda _{12}\,\mu _{11}\,\mu _{12}\,\mu _{21}\,\mu _{22}
+\lambda _{21}\,\lambda _{22}\,\mu _{11}\,\mu _{12}\,\mu _{21}\,\mu _{22}
\end{array}
\]
\[
E = \! \frac{A}{2} = (\mu _{11}\,\mu _{12}\,\mu _{22}+\mu _{11}\,\mu _{21}\,\mu _{22})(\lambda _{11}+\lambda _{12}+\lambda _{21}+\lambda _{22})
\]
\begin{multline*}
F = (\lambda _{11}+\lambda _{12}+\lambda _{21}+\lambda _{22})\,*\\
(\lambda _{21}\,\mu _{11}\,\mu _{12}\,\mu _{22}+\lambda _{22}\,\mu _{11}\,\mu _{12}\,\mu _{21}+\\
\lambda _{21}\,\mu _{11}\,\mu _{21}\,\mu _{22}+\mu _{11}\,\mu _{12}\,\mu _{21}\,\mu _{22})
\end{multline*}
\begin{multline*}
G = (\lambda _{11} \lambda _{12}+\lambda _{21}+\lambda _{22})*\\
\,(\lambda _{12}\,\mu _{11}\,\mu _{12}\,\mu _{22}+\lambda _{11}\,\mu _{12}\,\mu _{21}\,\mu _{22}+\\
\lambda _{12}\,\mu _{11}\,\mu _{21}\,\mu _{22}+\mu _{11}\,\mu _{12}\,\mu _{21}\,\mu _{22})
\end{multline*}
\begin{multline*}
H = (\lambda _{11}+\lambda _{12}+\lambda _{21}+\lambda _{22})\,*\\
(
\lambda _{12}\,\lambda _{21}\,\mu _{11}\,\mu _{12}\,\mu _{22}
+\lambda _{12}\,\lambda _{22}\,\mu _{11}\,\mu _{12}\,\mu _{21}\\
+\lambda _{11}\,\lambda _{21}\,\mu _{12}\,\mu _{21}\,\mu _{22} 
+\lambda _{12}\,\lambda _{21}\,\mu _{11}\,\mu _{21}\,\mu _{22}\\
+\lambda _{12}\,\mu _{11}\,\mu _{12}\,\mu _{21}\,\mu _{22}
+\lambda _{21}\,\mu _{11}\,\mu _{12}\,\mu _{21}\,\mu _{22} )
\end{multline*}

\section*{Appendix C: Optimal Policy Using Dynamic Programming for the $N=2$, $m=1$ Equally Weighted $w = 0.5$ System} \label{A3}

The $N=2$ $m=1$ equally weighed $w = 0.5$ dynamic linear program may be written as:
\begin{multline}
\max_{\bar{J},h(\mathcal{N}_1),h(\mathcal{N}_2),h(W_{1,1}),h(W_{1,2}),h(W_{2,1}),h(W_{2,2}),h(W_{1,1}^E),h(W_{2,1}^E)} \hspace{-9em}  \bar{J} \\
\mbox{s. t. } \bar{J} + h(\mathcal{N}_1) \le \frac{\sum_j \lambda_{2,j}}{\alpha} + \frac{\lambda_{1,1}}{\alpha}h(W_{1,1}) \\
+ \frac{\lambda_{1,2}}{\alpha}h(W_{1,2}) + \big ( 1 - \frac{\sum_j \lambda_{1,j}}{\alpha} \big ) h(\mathcal{N}_1 ) \\
 \bar{J} + h(\mathcal{N}_2) \le \frac{\sum_j \lambda_{1,j}}{\alpha} + \frac{\lambda_{2,1}}{\alpha}h(W_{2,1}) + \frac{\lambda_{2,2}}{\alpha}h(W_{2,2}) \\
 + \big ( 1 - \frac{\sum_j \lambda_{2,j}}{\alpha} \big ) h(\mathcal{N}_2 ) \\
 \bar{J} + h(W_{1,1}) \le 1 + \frac{\mu_{1,1}}{\alpha}h(\mathcal{N}_1) + \big ( 1 - \frac{\mu_{1,1}}{\alpha} \big ) h(W_{1,1} ) \\
  \bar{J} + h(W_{1,1}) \le 1 + \frac{\mu_{1,1}}{\alpha}h(W_{1,2}^E) + \big ( 1 - \frac{\mu_{1,1}}{\alpha} \big ) h(W_{1,1} ) \\
   \bar{J} + h(W_{1,2}) \le 1 + \frac{\mu_{1,2}}{\alpha}h(\mathcal{N}_2) + \big ( 1 - \frac{\mu_{1,2}}{\alpha} \big ) h(W_{1,2} ) \\
  \bar{J} + h(W_{1,2}) \le 1 + \frac{\mu_{1,2}}{\alpha}h(W_{2,1}^E) + \big ( 1 - \frac{\mu_{1,2}}{\alpha} \big ) h(W_{1,2} ) \\
   \bar{J} + h(W_{2,1}) \le 1 + \frac{\mu_{2,1}}{\alpha}h(\mathcal{N}_1) + \big ( 1 - \frac{\mu_{2,1}}{\alpha} \big ) h(W_{2,1} ) \\
  \bar{J} + h(W_{2,1}) \le 1 + \frac{\mu_{2,1}}{\alpha}h(W_{1,2}^E) + \big ( 1 - \frac{\mu_{2,1}}{\alpha} \big ) h(W_{2,1} ) \\
     \bar{J} + h(W_{2,2}) \le 1 + \frac{\mu_{2,2}}{\alpha}h(\mathcal{N}_2) + \big ( 1 - \frac{\mu_{2,2}}{\alpha} \big ) h(W_{2,2} ) \\
  \bar{J} + h(W_{2,2}) \le 1 + \frac{\mu_{2,2}}{\alpha}h(W_{2,1}^E) + \big ( 1 - \frac{\mu_{2,2}}{\alpha} \big ) h(W_{2,2} ) \\
     \bar{J} + h(W_{1,2}^E) \le 2 + \frac{\mu_{1,2}}{\alpha}h(\mathcal{N}_2) + \big ( 1 - \frac{\mu_{1,2}}{\alpha} \big ) h(W_{1,2}^E ) \\
  \bar{J} + h(W_{1,2}^E) \le 2 + \frac{\mu_{1,2}}{\alpha}h(W_{2,1}^E) + \big ( 1 - \frac{\mu_{1,2}}{\alpha} \big ) h(W_{1,2}^E ) \\
     \bar{J} + h(W_{2,1}^E) \le 2 + \frac{\mu_{2,1}}{\alpha}h(\mathcal{N}_1) + \big ( 1 - \frac{\mu_{2,1}}{\alpha} \big ) h(W_{2,1}^E ) \\
  \bar{J} + h(W_{2,1}^E) \le 2 + \frac{\mu_{2,1}}{\alpha}h(W_{1,2}^E) + \big ( 1 - \frac{\mu_{2,1}}{\alpha} \big ) h(W_{2,1}^E ) \\
\end{multline}
with uniformized rate $\alpha \ge \sum_{i} \sum_j \lambda_{i,j} + \max \{\mu_{i,j} \} $ from Figure \ref{chain}. 

The first two constraints deal with the idle AT states of $\mathcal{N}_1$ and $\mathcal{N}_2$ which by this definition have no controls but have an expected cost associated with the ratio of users they expect to miss i.e. the \% which arrive at the other region. Constraints 3-10 are for the en routes states $W_{1,1},W_{1,2},W_{2,1},$ and $W_{2,2}$ each with 2 controls: do nothing and let the AT become idle or send the empty AT to the other region with an associated cost of 1 because 100\% of expected user arrivals will be missed because the 1 and only AT is busy. Constraints 11-14 deal with states $W_{1,2}^E$ and $W_{2,1}^E$ which echo states $W_{1,2}$ and $W_{2,1}$  in terms of controls but have an associated cost of 2 because 100\% of users arriving will not obtain an AT and 100\% of time in that state is an AT driving empty.

\end{document}